\documentclass[12pt,a4paper]{article}

\title{The Injective Spectrum of a Right Noetherian Ring II:\\Sheaves and Torsion Theories}
\author{Harry Gulliver}
\date{}

\usepackage[
top    = 2.5cm,
bottom = 2.5cm,
left   = 2.3cm,
right  = 2.3cm]{geometry}

\usepackage{upgreek}
\let\tau\uptau

\usepackage{amssymb}
\usepackage{amsmath}
\usepackage{amsthm}

\usepackage{tikz}
\usetikzlibrary{arrows}

\usepackage{fancyhdr}
\usepackage{enumerate}

\renewcommand{\qed}{\hfill$\blacksquare$\gap}
\newcommand{\gap}{\par \vspace{5mm}}
\let\ideal\trianglelefteq

\let\take\smallsetminus

\newcommand{\Mod}{\mathrm{Mod}}
\mathchardef\mhyphen="2D

\newcommand{\tr}{\mathrm{Tr}}

\newcommand{\ann}{\mathrm{ann}}

\newcommand{\spec}{\mathrm{Spec}}
\renewcommand{\mod}{\mathrm{mod}}

\newcommand{\F}{\mathcal{F}}
\newcommand{\T}{\mathcal{T}}

\let\take\smallsetminus

\newcommand{\injspec}{\mathrm{InjSpec}}
\newcommand{\cl}{\mathrm{cl}}

\newcommand{\torspec}{\mathrm{TorSpec}}

\newcommand{\fp}{\mathrm{fp}}

\newtheorem{thm}{Theorem}[section]
\renewcommand{\proof}[1][]{\noindent\textsc{Proof:} \textit{#1}\par}

\newtheorem{lemma}[thm]{Lemma}
\newtheorem{cor}[thm]{Corollary}
\newtheorem{ex}[thm]{Example}
\newtheorem{prop}[thm]{Proposition}
\newtheorem{qn}{Question}

\newcommand{\sh}{\mathrm{Sh}}
\newcommand{\qcoh}{\mathrm{QCoh}}
\newcommand{\res}{\mathrm{res}}
\newcommand{\tors}{\mathrm{tors}}

\newcommand{\uline}[1]{\underline{\smash{#1}}}

\begin{document}

\maketitle
\pagenumbering{arabic}

\begin{abstract}
This is the second of two papers on the injective spectrum of a right noetherian ring. In \cite{gulliver1}, we defined the injective spectrum as a topological space associated to a ring (or, more generally, a Grothendieck category), which generalises the Zariski spectrum. We established some results about the topology and its links with Krull dimension, and computed a number of examples.

In the present paper, which can largely be read independently of the first, we extend these results by defining a sheaf of rings on the injective spectrum and considering sheaves of modules over this structure sheaf and their relation to modules over the original ring. We then explore links with the spectrum of prime torsion theories developed by Golan \cite{golan} and use this torsion-theoretic viewpoint to prove further results about the topology.
\end{abstract}

\tableofcontents

\section{Introduction and Background}

\subsection{Conventions}

Throughout, all rings will be associative and unital, but not necessarily commutative, and all modules will be unital right modules, unless otherwise specified. If $R$ is a ring, we denote by $\Mod\mhyphen R$ the category of all right $R$-modules, and by $\mod\mhyphen R$ the full subcategory of finitely presented modules. If $M$ is a module, or more generally an object in a Grothendieck category, we denote by $E(M)$ an injective hull of $M$. For modules (objects of a Grothendieck category) $L$ and $M$, we denote by $(L,M)$ the group of maps $L\to M$.

By ``functor'' we always mean ``additive, covariant functor''. For a Grothendieck category $\mathcal{A}$, we denote by $\mathcal{A}^\fp$ the full subcategory of finitely presented objects; so $\mod\mhyphen R = (\Mod\mhyphen R)^\fp$.

\subsection{Torsion Theories}

Recall (\textit{e.g.}, \cite[\S11.1.1]{prest}) that a {\bf torsion theory} in a Grothendieck category $\mathcal{A}$ is a pair $(\T,\F)$ of classes of objects such that there are no non-zero maps from objects of $\T$ to objects of $\F$, and both classes are maximal with respect to this property. This is equivalent to $\T$ being closed under quotients, extensions, and arbitrary coproducts, and there being no maps from $\T$ to $\F$, equivalently to $\F$ being closed under subobjects, extensions, and arbitrary products, and there being no maps from $\T$ to $\F$. In a torsion theory, $\T$ is called the torsion class, and $\F$ the torsionfree class.

A torsion theory is {\bf hereditary} if $\T$ is closed under subobjects, equivalently if $\F$ is closed under injective hulls. We shall consider here only hereditary torsion theories, and so shall henceforth omit the adjective ``hereditary''. A {\bf Serre subcategory} in an abelian category is a full subcategory which is closed under subobjects, quotients, and extensions; so a hereditary torsion class is precisely a Serre subcategory which is additionally closed under coproducts.

An alternative description of a hereditary torsion theory on $\mathcal{A}$ is given by the {\bf torsion radical} or {\bf torsion functor}. This is the subfunctor $\tau$ of the identity functor on $\mathcal{A}$ such that for any object $A$, $\tau(A)$ is the largest subobject of $A$ contained in $\T$. Conversely, given a left exact subfunctor $\tau$ of the identity functor such that $\tau(A/\tau(A))=0$ for all objects $A$, then setting
\[\T=\{A\in\mathcal{A}\mid \tau(A)=A\},\quad \F=\{A\in\mathcal{A}\mid \tau(A)=0\}\]
gives a torsion theory; see \cite[Chapter VI]{stenstrom} for more details.

The principal significance of Serre subcategories and torsion theories comes from the following:

\begin{prop}[See Chapter 4, especially Sections 4.3 and 4.4, of \cite{popescu}]\label{localisation}
Let $\mathcal{A}$ be an abelian category and $\mathcal{S}$ a Serre subcategory. Then:
\begin{enumerate}
\item There exist an abelian category $\mathcal{A}/\mathcal{S}$ and a dense, exact functor $Q_\mathcal{S}:\mathcal{A}\to\mathcal{A}/\mathcal{S}$ with kernel $\mathcal{S}$ obeying the following universal property:

\begin{center}
\begin{tikzpicture}
\node (A) at (0,1.5) {$\mathcal{A}$};
\node (AS) at (0,0) {$\mathcal{A}/\mathcal{S}$};
\node (B) at (4,0) {$\mathcal{B}$};
\draw
(A) edge[->,>=angle 90] node[left] {$Q_\mathcal{S}$} (AS)
(A) edge[->,>=angle 90] node[above right] {$F$} (B)
(AS) edge[dashed,->,>=angle 90] node[below] {$\hat{F}$} (B);
\end{tikzpicture}
\end{center}
whenever $\mathcal{B}$ is an abelian category and $F:\mathcal{A}\to\mathcal{B}$ is an exact functor such that $F(A)=0$ for all $A\in\mathcal{S}$, then there exists a unique exact functor $\hat{F}:\mathcal{A}/\mathcal{S}\to \mathcal{B}$ such that $F=\hat{F}\circ Q_\mathcal{T}$.

\item If $\mathcal{A}$ is Grothendieck, then $\mathcal{S}$ is closed under coproducts (\textit{i.e.}, is a torsion class in $\mathcal{A}$) if and only if $Q_\mathcal{S}$ admits a right adjoint, which we denote $i_\mathcal{S}$.

\item If $\mathcal{A}$ is Grothendieck and $\mathcal{S}$ is a torsion class, then $i_\mathcal{S}$ is fully faithful and $\mathcal{A}/\mathcal{S}$ is Grothendieck. Moreover, for any object $A\in\mathcal{A}$, the localisation $i_\mathcal{S}Q_\mathcal{S}(A)$ can be described as
\[i_\mathcal{S}Q_\mathcal{S}(A)=\pi^{-1}\left(\tau_\mathcal{S}\left(\frac{E(A/\tau_{\mathcal{S}}(A))}{A/\tau_\mathcal{S}(A)}\right)\right),\]
where
\[\pi:E(A/\tau_\mathcal{S}(A))\to \frac{E(A/\tau_\mathcal{S}(A))}{A/\tau_\mathcal{S}(A)}\]
is the quotient map. That is, to localise an object of a Grothendieck category at a torsion class, we quotient out the torsion part to obtain a torsionfree object, then look at the part of the injective hull which becomes torsion modulo this torsionfree object.
\end{enumerate}
\end{prop}

We call $\mathcal{A}/\mathcal{S}$ the {\bf quotient category} or {\bf localisation} of $\mathcal{A}$ by $\mathcal{S}$, $Q_\mathcal{S}$ the {\bf quotient functor} or {\bf localisation functor}, and $i_\mathcal{S}$ the {\bf adjoint inclusion functor}.\gap

A torsion theory is of {\bf finite type} if it satisfies the equivalent conditions of the following

\begin{lemma}[\cite{prest}, 11.1.12, 11.1.14, 11.1.26]\label{finite type}
Let $(\mathcal{T},\mathcal{F})$ be a torsion theory in the Grothen\-dieck category $\mathcal{A}$. Then the following are equivalent:
\begin{enumerate}
\item $\tau_\mathcal{T}$ commutes with directed colimits.
\item $i_\mathcal{T}$ commutes with directed colimits of monomorphisms.
\item $\mathcal{F}$ is closed under directed colimits.
\end{enumerate}
Moreover, if $\mathcal{A}$ is locally finitely presented (\textit{i.e.}, has a generating set of finitely presented objects), then these conditions are equivalent to $\mathcal{T}$ being generated as a torsion class by $\mathcal{T}\cap\mathcal{A}^\mathrm{fp}$, the finitely presented torsion objects, and this establishes a bijection between Serre subcategories of $\mathcal{A}^\fp$ and torsion classes of finite type in $\mathcal{A}$.

When these equivalent conditions hold, if $A\in\mathcal{A}$ is finitely generated, then $Q_\mathcal{T}(A)\in\mathcal{A}/\mathcal{T}$ is finitely generated, and if $\mathcal{G}$ is a generating family for $\mathcal{A}$, then $Q_\mathcal{T}\mathcal{G}$ is a generating family for $\mathcal{A}/\mathcal{T}$. If $\mathcal{A}$ is locally finitely presented, then ``finitely generated'' can be replaced by ``finitely presented'' in this paragraph.
\end{lemma}

It is easy to see from the closure conditions that any intersection of torsion (resp. torsionfree) classes is itself a torsion (resp. torsionfree) class. Therefore, given an indexing set $I$ and a torsion theory $(\mathcal{T}_i,\mathcal{F}_i)$ for each $i\in I$, we can construct two new torsion theories. The first of these has torsion class $\bigcap_{i\in I}\mathcal{T}_i$; we denote the torsionfree class for this theory $\sum_{i\in I}\mathcal{F}_i$. The second has torsionfree class $\bigcap_{i\in I}\mathcal{F}_i$; we denote its torsion class $\sum_{i\in I}\mathcal{T}_i$.

It is not hard to check that these intersection and sum operations make the set of torsion classes (partially ordered under inclusion), into a complete lattice. Similarly, the set of torsionfree classes is a complete lattice, and these two lattices are dual to each other. See \cite[\S 1]{golan} for details, though be aware that the notation there differs significantly from here.

Given a class $C$ of objects in $\mathcal{A}$, we denote by $\mathcal{T}(C)$ the intersection of all torsion classes containing $C$ and call it the {\bf torsion class generated by $C$}. Similarly, we denote by $\mathcal{F}(C)$ the intersection of all torsionfree classes containing $C$ and call it the {\bf torsionfree class cogenerated by $C$}. When $C=\{A\}$ consists of a single object, we omit the braces, writing simply $\mathcal{T}(A)$ and $\mathcal{F}(A)$.

The following useful result is well known.

\begin{lemma}\label{cogenerating}
Let $\mathcal{A}$ be a locally noetherian Grothendieck category and $S$ a set of objects of $\mathcal{A}$. Let $E(S)$ denote the set of injective hulls of objects in $S$, and $P$ the product of all objects in $S$. Then $\mathcal{F}(S)$ consists of all subobjects of direct products of objects of $E(S)$, $\mathcal{F}(S)=\mathcal{F}(E(S))=\mathcal{F}(P)$, and $\mathcal{T}_{\mathcal{F}(S)}$ consists of those objects $A$ such that $(A,E)=0$ for all $E\in E(S)$.
\end{lemma}

\subsection{The Injective Spectrum; Prior Results}\label{injspec intro}

Beyond the basic definitions, this paper is largely independent of the preceding paper \cite{gulliver1}. We recall here the relevant definitions and a small number of results from \cite{gulliver1}, which will be relevant to this paper, and may be viewed as ``black box'' results for the reader who is more interested in this paper alone than in \cite{gulliver1}.

Let $\mathcal{A}$ be a Grothendieck abelian category. The {\bf injective spectrum} of $\mathcal{A}$, denoted $\injspec(\mathcal{A})$, is the set of isoclasses of indecomposable injective objects of $\mathcal{A}$, topologised as follows. For any finitely presented object $A\in\mathcal{A}$, write $[A]$ for the set of indecomposable injectives $E$ such that $(A,E)=0$; take the set of all $[A]$ as $A$ ranges over $\mathcal{A}^\fp$ as a basis of open sets for a topology on $\injspec(\mathcal{A})$, which we call the Zariski topology.

For $A\in\mathcal{A}^\fp$, write $(A)$ for the set of indecomposable injectives $E$ such that $(A,E)\neq 0$; \textit{i.e.}, the complement of $[A]$ in $\injspec(\mathcal{A})$. We refer to such sets as {\bf basic closed} sets for the Zariski topology on $\injspec(\mathcal{A})$. If $\mathcal{A}$ is locally noetherian (\textit{i.e.}, has a generating set of noetherian objects), then there is an alternative topology on $\injspec(\mathcal{A})$, called the {\bf Ziegler topology}, having the sets $(A)$ for $A\in\mathcal{A}^\fp$ as a basis of open sets. The sets $(A)$ for $A\in\mathcal{A}^\fp$ are precisely the compact open sets of the Ziegler topology on $\injspec(\mathcal{A})$. See \cite[\S\S 5.6, 14.1]{prest} for more details.

When we refer to the injective spectrum without specifying a topology, we shall always mean the Zariski topology; however, on occasion we shall find it useful to switch to the Ziegler topology in proofs.

In the event that $\mathcal{A}=\Mod\mhyphen R$ is the category of right modules over some ring $R$, we write simply $\injspec(R)$ as shorthand for $\injspec(\Mod\mhyphen R)$. We have the following result of Gabriel, who first considered the injective spectrum.

\begin{thm}[\cite{gabriel}, \S VI.3]
Let $R$ be a commutative noetherian ring. Then there is a homeomorphism $\injspec(R)\cong\spec(R).$
\end{thm}

If $E,F\in\injspec(\mathcal{A})$ are indecomposable injectives, we write $E\leadsto F$ and say that $E$ {\bf specialises to} $F$ if $F\in\cl(E)$; \textit{i.e.}, if every closed set containing $E$ also contains $F$. The following is an adaptation of a result from \cite{gulliver1}.

\begin{lemma}[\cite{gulliver1}, Lemma 2.1]\label{specialisation tfae}
Let $\mathcal{A}$ be a locally noetherian Grothendieck category. For $E,F\in \injspec(\mathcal{A})$, the following are equivalent:
\begin{enumerate}
\item $E\leadsto F$;
\item $E\in\mathcal{F}(F)$;
\item $\mathcal{F}(E)\subseteq\mathcal{F}(F)$.
\end{enumerate}
\end{lemma}

\proof
The equivalence between (1) and (2) comes from parts (1) and (4) of \cite[Lemma 2.1]{gulliver1}, rephrased in terms of torsionfree classes. The equivalence between (2) and (3) is by definition of $\mathcal{F}(E)$.\qed

We shall also require the following results about the topology of the injective spectrum.

\begin{prop}[\cite{gulliver1}, Corollary 3.5]\label{t0}
For any locally noetherian Grothendieck category $\mathcal{A}$, $\injspec(\mathcal{A})$ is $T_0$, \textit{i.e.}, Kolmogorov.
\end{prop}

\begin{lemma}[\cite{gulliver1}, Lemma 3.10]\label{direct sum basic open sets}
Let $\mathcal{A}$ be any Grothendieck category and
\[0\to A\to C\to B\to 0\]
a short exact sequence of finitely presented objects in $\mathcal{A}$. Then $(C)=(A)\cup (B)$.
\end{lemma}

\begin{thm}[\cite{gulliver1}, Theorem 3.15]\label{domain has irred spec}
If $R$ is a right noetherian domain, then $\injspec(R)$ is irreducible and has $E(R_R)$ as a generic point.
\end{thm}

Finally, a technical Lemma which does not appear explicitly in \cite{gulliver1}, but follows from results there concerning Krull dimension and critical dimension.

\begin{lemma}\label{crit subobj}
Let $\mathcal{A}$ be a Grothendieck category and $A$ a non-zero noetherian object in $\mathcal{A}$. Then there is a non-zero subobject $B$ of $A$ with the property that for any proper quotient $B/C$, there are no non-zero morphisms $B/C\to E(B)$.
\end{lemma}

\proof
Being noetherian, $A$ has a critical subobject $B$ (in the sense of Krull dimension) \cite[\S 6.2]{mcr}. Any such $B$ has the required property. For given any proper quotient $B/C$, we have $K(B/C)<K(B)$, by definition, and so $(B/C,E(B))=0$ by \cite[3.1.4 \& 3.2.6]{gulliver1}.\qed

\subsection{Outline of Paper}

This paper and its predecessor \cite{gulliver1} together present the results of the author's PhD thesis, which was prepared under the supervision of Prof Mike Prest and submitted to the University of Manchester in June 2019. This paper is written to be independent of the prequel, so while a knowledge of that may be helpful in understanding parts of the present paper, it is not necessary.

In section \ref{sheaves}, we consider a sheaf of rings on the injective spectrum of a ring, originally constructed by Gabriel. We show that the ring of global sections of this sheaf is not always isomorphic (or even Morita equivalent) to the original ring, but that it is if the ring is a noetherian domain. We then construct two functors from $R$-modules to sheaves of modules over $\injspec(R)$, and establish a necessary and sufficient criterion for when these functors coincide, as well as proving that when they do the resulting sheaves are quasicoherent.

In section \ref{torsion}, we consider an alternative topological space, the torsion spectrum, introduced by Golan. We show that this is homeomorphic to the injective spectrum in any locally noetherian Grothendieck category, and then exploit this connection to prove further results about torsion theories and sobriety of the injective spectrum in its Ziegler topology.

Finally, in section \ref{spectral and noetherian}, we consider whether the injective spectrum is a spectral space, and show that if it is, we can isolate basic closed sets as spectra of related Grothendieck categories. If the injective spectrum is also noetherian, then this extends to all closed sets.

\subsection{Acknowledgements}

I owe an enormous debt of gratitude to Prof Mike Prest, who supervised my PhD, in which I completed the work which forms this paper; his guidance and patience have been exemplary. Thanks are also due to Tommy Kucera, Ryo Kanda, Lorna Gregory, and Marcus Tressl, as well as my examiners, Omar Le{\'o}n S{\'a}nchez and Gwyn Bellamy, for stimulating discussions and helpful feedback. Finally, I am grateful to EPSRC for providing the funding for my PhD.

\section{Sheaves}\label{sheaves}

\subsection{The Structure Sheaf}

We describe a sheaf of rings on $\injspec(R)$, which was developed along with the topology by Gabriel in \cite[\S VI.3]{gabriel}, though our presentation is rather different. We begin by constructing a presheaf-on-a-basis. Given a basic open set $[M]$ (for some $M\in\mod\mhyphen R$), associate the torsion class $\T(M)$. This depends only on the set $[M]$, not on the choice of representing module $M$; for it follows from Lemma \ref{cogenerating} that the associated torsionfree class is $\F([M])$, cogenerated by the indecomposable injectives in $[M]$, so if $[M]=[N]$, then $\F([M])=\F([N])$, and so $\T(M)=\T(N)$. Since $M$ is finitely presented, $\T(M)$ is of finite type.

For convenience, we denote the localised category $(\Mod\mhyphen R)/\T(M)$ by $(\Mod\mhyphen R)_M$ and the localisation functor by $Q_M$. Let $R_M$ denote the endomorphism ring of $Q_M(R_R)$. This will be the ring associated to the basic open set $[M]$ by our presheaf-on-a-basis. If $[N]\subseteq[M]$, then $\mathcal{F}_N$ is cogenerated by a subset of (a cogenerating set of) $\mathcal{F}_M$, so $\mathcal{F}_N\subseteq\mathcal{F}_M$. Therefore $\mathcal{T}(M)\subseteq\mathcal{T}(N)$; so, by the universal property of localisation (Proposition \ref{localisation}), the quotient functor $Q_N$ factors through $Q_M$ by a unique exact functor $Q_{M,N}$:

\begin{center}
\begin{tikzpicture}
\node (A) at (0,1.5) {$\Mod\mhyphen R$};
\node (M) at (0,0) {$(\Mod\mhyphen R)_M$};
\node (N) at (4,0) {$(\Mod\mhyphen R)_N$};
\draw
(A) edge[->,>=angle 90] node[left] {$Q_M$} (M)
(A) edge[->,>=angle 90] node[above right] {$Q_N$} (N)
(M) edge[->,>=angle 90] node[below] {$Q_{M,N}$} (N);
\end{tikzpicture}
\end{center}

Therefore $Q_{M,N}Q_MR=Q_NR$, and so $R_N=(Q_{M,N}Q_MR,Q_{M,N}Q_MR)$. Moreover, $Q_{M,N}$ gives a ring map $R_M=(Q_MR,Q_M,R)\to (Q_{M,N}Q_MR,Q_{M,N}Q_MR)=R_N$, which we denote $\rho_{M,N}$. So for each inclusion of basic open sets $[N]\subseteq [M]$, we have a restriction map $\rho_{M,N}:R_M\to R_N$.

Similarly, if $[L]\subseteq[N]\subseteq[M]$, $Q_{M,L}:(\Mod\mhyphen R)_M\to (\Mod\mhyphen R)_L$ is the unique exact functor such that this diagram commutes:

\begin{center}
\begin{tikzpicture}
\node (A) at (0,1.5) {$\Mod\mhyphen R$};
\node (M) at (0,0) {$(\Mod\mhyphen R)_M$};
\node (L) at (4,0) {$(\Mod\mhyphen R)_L$};
\draw
(A) edge[->,>=angle 90] node[left] {$Q_M$} (M)
(A) edge[->,>=angle 90] node[above right] {$Q_L$} (L)
(M) edge[->,>=angle 90] node[below] {$Q_{M,L}$} (L);
\end{tikzpicture}
\end{center}

But the following diagram commutes, and the bottom row is exact:

\begin{center}
\begin{tikzpicture}
\node (A) at (0,2) {$\Mod\mhyphen R$};
\node (M) at (0,0) {$(\Mod\mhyphen R)_M$};
\node (N) at (4,0) {$(\Mod\mhyphen R)_N$};
\node (L) at (8,0) {$(\Mod\mhyphen R)_L$};
\draw
(A) edge[->,>=angle 90] node[left] {$Q_M$} (M)
(A) edge[->,>=angle 90] node[below left] {$Q_N$} (N)
(M) edge[->,>=angle 90] node[below] {$Q_{M,N}$} (N)
(N) edge[->,>=angle 90] node[below] {$Q_{N,L}$} (L)
(A) edge[->,>=angle 90] node[above right] {$Q_L$} (L);
\end{tikzpicture}
\end{center}

Therefore $Q_{M,L}=Q_{N,L}\circ Q_{M,N}$. In particular, $\rho_{M,L}=\rho_{N,L}\circ\rho_{M,N}$. Therefore the assignment taking a basic open set $[M]$ to $R_M$ and an inclusion of basic open sets $[N]\subseteq[M]$ to $\rho_{M,N}$ is a presheaf-on-a-basis on $\injspec(R)$. This is sufficient for the sheafification process to work \cite[\S3.2]{ega}, and so we obtain a sheaf of rings $\mathcal{O}_R$ on $\injspec(R)$, which we call the {\bf sheaf of finite type localisations}, or simply the {\bf structure sheaf}.

Of course, we must compare this to the usual structure sheaf in the commutative case. Indeed, we have the following:

\begin{thm}[\cite{gabriel}, \S VI.3]\label{gabriel sheaf = zariski sheaf}
If $R$ is commutative noetherian and $\injspec(R)$ is identified with $\spec(R)$ via the Matlis bijection, then the sheaf of finite type localisations is isomorphic to the usual Zariski structure sheaf.
\end{thm}

We have presented the sheaf of finite type localisations specifically over a ring, whereas we have defined the injective spectrum for an arbitrary Grothendieck category. The construction will still work for any locally noetherian Grothendieck category, by choosing a generator to stand in place of $R_R$; however, since there is not generally a canonical choice of generator in a Grothendieck category, this becomes non-canonical. As such, throughout this section we shall stick to the case of rings.\gap

A key property of the Zariski spectrum of a commutative ring is that the ring of global sections of the structure sheaf is simply the original ring. We begin with an example to show that this can fail for the injective spectrum.

\begin{ex}\label{kA2 counterexample}
Let $k$ be a field and $R=kA_2$ be the path algebra over $k$ of the quiver $A_2:1\to 2$; then the ring of global sections of the sheaf of finite type localisations is $k\oplus M_2(k)$, not $R$.
\end{ex}

By standard results on quiver representations, $R$ has exactly two indecomposable injectives, namely the representations $(k\to 0)$ and $(k\to k)$; the topology on the injective spectrum is discrete, by \cite[Proposition 4.1]{gulliver1} or an easy calculation. The ring of global sections is therefore simply the direct sum of the two stalks.

For the stalk at $(k\to 0)$, we take the torsionfree class cogenerated by $(k\to 0)$ and localise $R_R=(k\to k)\oplus (0\to k)$; obtaining $(k\to 0)$. We then take the endomorphism ring of this, which is simply $k$. For the stalk at $(k\to k)$, we localise $R_R$ at the torsionfree class cogenerated by $(k\to k)$, obtaining $(k\to k)^2$, which has endomorphism ring $M_2(k)$, the $2\times 2$ matrix ring.

So the ring of global sections over $\injspec(kA_2)$ is $k\oplus M_2(k)$, which is not isomorphic - or even Morita equivalent - to $kA_2$.\qed

Having shown that, even for a very straightforward ring, the structure sheaf on the injective spectrum can fail to fulfill our expectations from the commutative case, we now show that it nonetheless often does.

\begin{thm}
Let $R$ be a right noetherian domain. Then the ring of global sections of the structure sheaf of $\injspec(R)$ is precisely $R$.
\end{thm}

\proof
By Theorem \ref{domain has irred spec}, $\injspec(R)$ is irreducible and has $E(R_R)$ as a generic point. It follows from Lemma \ref{specialisation tfae} that $E(R_R)$, and hence $R_R$, are torsionfree for every non-trivial torsion theory. Therefore, applying the localisation formula of Proposition \ref{localisation}, the localisation of $R$ at any torsion theory is the largest submodule of $E(R)$ which becomes torsion modulo $R$. So the presheaf-on-a-basis of localisations of $R$ associates to each basic open set a submodule of $E(R)$ (with the structure of a ring via its endomorphisms), and the restriction maps are simply inclusions into ever larger submodules of $E(R)$.

So for any global section $\sigma$ of the structure sheaf, and any point $E$, there is an open set $U_E\ni E$ and an element $e_E\in E(R)$ such that $\sigma(F)=e_E$ for all $F\in U_E$. But $\injspec(R)$ is irreducible, so given any two points $E$ and $F$, $U_E\cap U_F\neq \varnothing$; therefore there exists some $G\in U_E\cap U_F$, and $\sigma(E)=\sigma(G)=\sigma(F)$. So in fact there is a single element $e\in E(R)$ such that $\sigma(E)=e$ for all $E\in \injspec(R)$.

It remains to show that $e\in R$. For any $E\in \injspec(R)$, the submodule of $E(R)/R$ generated by $e+R$ must be $\mathcal{F}(E)$-torsion, by the localisation formula, so $((e+R)R,E)=0$ for all $E\in\injspec(R)$. Let $I=\ann_R(e+R)$, so $(e+R)R\cong R/I$. If $e\notin R$, then $I\neq R$, so $R/I$ is non-zero; but then $R/I$ has a simple quotient $S$, so $(R/I,E(S))\neq 0$, a contradiction. So indeed $e\in R$.\qed

\subsection{Sheaves of Modules}

Let $\mathcal{O}_R$ denote the sheaf of finite-type localisations on $\injspec(R)$. We now consider sheaves of $\mathcal{O}_R$-modules; let $\sh(R)$ denote the category of all such sheaves, with morphisms of sheaves as arrows.

We begin by describing two functors $\Mod\mhyphen R\to\sh(R)$. The first we call the {\bf tensor sheaf functor}; for $M\in\Mod\mhyphen R$, we write $\mathcal{M}_\otimes$ for the tensor sheaf of $M$, and for $f:M\to N$ a map of $R$-modules, we write $f_\otimes:\mathcal{M}_\otimes\to\mathcal{N}_\otimes$ for the induced map.

The tensor sheaf functor is defined as follows: given an $R$-module $M$ and an open set $U$ in $\injspec(R)$, we form $M\otimes_R \mathcal{O}_R(U)$. Given $U\subseteq V$ an inclusion of open sets in $\injspec(R)$, the restriction map $\mathcal{O}_R(V)\to\mathcal{O}_R(U)$ induces a restriction map $M\otimes_R \mathcal{O}_R(V)\to M\otimes_R\mathcal{O}_R(U)$. This gives a presheaf of $\mathcal{O}_R$-modules, whose sheafification we define to be $\mathcal{M}_\otimes$, the tensor sheaf of $M$.

Given $f:M\to N$ in $\Mod\mhyphen R$, we have for each open set $U$ an induced map $f\otimes_R\mathcal{O}_R(U):M\otimes_R\mathcal{O}_R(U)\to N\otimes_R\mathcal{O}_R(U)$. Since this acts on the first factor of the tensor product and the restriction maps act on the second factor, these two maps commute, and so we have a morphism of presheaves. Sheafification then gives a morphism of sheaves $f_\otimes:\mathcal{M}_\otimes\to\mathcal{N}_\otimes$.

The fact that this tensor sheaf construction is functorial is trivial to verify.

It will be useful at times to reach this functor by a slightly different route. Since $\mathcal{O}_R$ is the sheafification of the presheaf-on-a-basis $[A]\mapsto R_A$, we can form a presheaf-on-a-basis of modules $[A]\mapsto M\otimes_R R_A$ for any $M\in\Mod\mhyphen R$, and then sheafify this and extend to the whole topology to obtain a sheaf of $\mathcal{O}_R$-modules.

To see that these two constructions of the sheaf $\mathcal{M}_\otimes$ are the same, we consider the presheaf-on-a-basis $[A]\mapsto M\otimes_R R_{A}$ and the sheaf-on-a-basis $[A]\mapsto M\otimes_R \mathcal{O}_R([A])$. Since $\mathcal{O}_R$ is the sheafification of $[A]\mapsto R_{A}$, there is a natural map from $[A]\mapsto R_{A}$ to $\mathcal{O}_R$, which becomes an isomorphism when sheafified. Tensoring with $M$ gives a natural map from the presheaf $[A]\mapsto M\otimes_R R_A$ to the presheaf $[A]\mapsto M\otimes_R \mathcal{O}_R([A])$, which becomes an isomorphism when sheafified. Therefore the sheaf-on-a-basis version of these two sheaves associated to $M$ are canonically isomorphic, and hence so too are the full sheaves.\gap

Our second functor $\Mod\mhyphen R\to\sh(R)$ we call the {\bf torsion sheaf functor}; for $M\in\Mod\mhyphen R$, we write $\mathcal{M}_\mathrm{tors}$ for the torsion sheaf of $M$, and for $f:M\to N$ a map of $R$-modules, we write $f_\tors:\mathcal{M}_\tors\to\mathcal{N}_\tors$ for the induced map.

To define the torsion sheaf functor, we first consider torsion-theoretic localisation of modules. Recall the construction of the sheaf of finite type localisations $\mathcal{O}_R$. Given a torsion theory $\T$ in $\Mod\mhyphen R$ (particularly one of the form $\T(M)$ for $M\in\mod\mhyphen R$), we can take a ring $R_\T=(Q_\T R,Q_\T R)$, the endomorphism ring of the image of $R_R$ in the quotient category. There is a natural isomorphism of abelian groups $R_\T\cong i_\T Q_\T R$, and a canonical ring map $R\to R_\T$. The sheaf of finite type localisations was obtained by sheafifying the presheaf-on-a-basis $[A]\mapsto R_A$.

We will now mirror this construction with modules to obtain an $R_\T$-module structure on $M_\T:=(Q_\T R, Q_\T M)$ for any $M\in\Mod\mhyphen R$. This will gives us a presheaf-on-a-basis of modules over the presheaf-on-a-basis of rings of finite type localisations, which is enough information for sheafification to give us a sheaf of modules over $\mathcal{O}_R$.

\begin{lemma}
Let $\T$ be a torsion class in $\Mod\mhyphen R$. Then there is a functor $(-)_\T:\Mod\mhyphen R\to \Mod\mhyphen R_\T$. Moreover, if $\mathcal{S}\subseteq \T$ is an inclusion of torsion classes, then there are restriction maps of abelian groups $\res^M_{\mathcal{S},\T}:M_\mathcal{S}\to M_\T$ such that $\res^R_{\mathcal{S},\T}$ is a ring map and each $\res^M_{\mathcal{S},\T}$ is $R_\mathcal{S}$-linear when $M_\T$ has its $R_\mathcal{S}$ structure given by the ring map $\res^R_{\mathcal{S},\T}$.
\end{lemma}

\proof
We set $M_\T=(Q_\T R,Q_\T M)$. Since $R_\T=(Q_\T R,Q_\T R)$, we have a pairing $M_\T\times R_\T\to M_\T:(\mu,\rho)\to\mu\circ\rho$. This satsifies the axioms to make $M_\T$ into an $R_\T$-module by preadditivity of $(\Mod\mhyphen R)/\T$.

Given $f:M\to N$ in $\Mod\mhyphen R$, define $f_\T:M_\T\to N_\T:\mu\mapsto (Q_\T f)\circ\mu$. It is trivial from functoriality of $Q_\T$ and preadditivity to verify that this defines an additive functor $\Mod\mhyphen R\to\Mod\mhyphen R_\T$.

Now let $\mathcal{S}\subseteq\T$ be an inclusion of torsion classes. By the universal property of torsion-theoretic localisation (Proposition \ref{localisation}), there is a unique exact functor $Q_{\mathcal{S},\T}:(\Mod\mhyphen R)/\mathcal{S}\to(\Mod\mhyphen R)/\T$ such that $Q_\T=Q_{\mathcal{S},\T}\circ Q_\mathcal{S}$. Therefore $Q_{\mathcal{S},\T}$ induces maps of abelian groups $(Q_\mathcal{S} R,Q_\mathcal{S} R)\to (Q_\T R, Q_\T R)$ and $(Q_\mathcal{S} R,Q_\mathcal{S} M)\to (Q_\T R, Q_\T M)$; \textit{i.e.}, $R_\mathcal{S}\to R_\T$ and $M_\mathcal{S}\to M_\T$. These are our restriction maps. Again, the linearity properties follow easily from properties of the functors.\qed

Now, given an $R$-module $M$, we define a presheaf-on-a-basis by assigning to the basic open set $[A]$ the $R_{A}$-module $M_{\T(A)}$. By the above Lemma, this is indeed a presheaf. Its sheafification is the torsion sheaf $\mathcal{M}_\tors$ associated to $M$.

Given a map $f:M\to N$ in $\Mod\mhyphen R$, we construct a morphism of presheaves between the presheaves-on-a-basis associated to $M$ and $N$. Sheafification then turns this into a morphism of sheaves. For each basic open set $[A]$, we have $f_{\T(A)}:M_{\T(A)}\to N_{\T(A)}$ as in the Lemma. We need to check that these cohere with the restriction maps; \textit{i.e.}, that if $[A]\subseteq [B]$ (so $\T(B)\subseteq \T(A))$, we have a commuting diagram

\begin{center}
\begin{tikzpicture}

\node (mb) at (0,2) {$M_{\T(B)}$};
\node (ma) at (0,0) {$M_{\T(A)}$};
\node (nb) at (3,2) {$N_{\T(B)}$};
\node (na) at (3,0) {$N_{\T(A)}$};

\draw
(mb) edge[->,>=angle 90] node[left] {$\res^M_{\T(B),\T(A)}$} (ma)
(nb) edge[->,>=angle 90] node[right] {$\res^N_{\T(B),\T(A)}$} (na)
(mb) edge[->,>=angle 90] node[above] {$f_{\T(B)}$} (nb)
(ma) edge[->,>=angle 90] node[below] {$f_{\T(A)}$} (na);

\end{tikzpicture}
\end{center}

To show that this diagram does indeed commute, take $\mu\in M_{\T(B)}$. Following the diagram anticlockwise, $\mu$ maps first to $Q_{\T(B),\T(A)}(\mu)\in M_{\T(A)}$, then to $(Q_{\T(A)} f)\circ (Q_{\T(B),\T(A)} \mu)\in N_{\T(A)}$. Following clockwise, $\mu$ maps first to $(Q_{\T(B)} f)\circ \mu\in N_{\T(B)}$, which then maps to $Q_{\T(B),\T(A)}((Q_{\T(B)} f)\circ \mu)=(Q_{\T(A)} f)\circ (Q_{\T(B),\T(A)}\mu)$, the same as when going anticlockwise.

So we have a functor from $\Mod\mhyphen R$ to presheaves-on-a-basis of modules. Sheafifying then gives the desired torsion sheaf functor $\Mod\mhyphen R\to \sh(R):M\mapsto \mathcal{M}_\tors$.\gap

So we have two functors $\Mod\mhyphen R\to \sh(R)$, the torsion sheaf functor and the tensor sheaf functor. We now consider the relationship between them.

\begin{lemma}
For any $M$ in $\Mod\mhyphen R$ and $\T$ a torsion class in $\Mod\mhyphen R$, there is a morphism of $R_\T$-modules $\theta_{M,\T}:M\otimes_R R_\T\to M_\T$. If $\mathcal{S}\subseteq\mathcal{T}$, then there is a commuting diagram

\begin{center}
\begin{tikzpicture}

\node (tens s) at (0,2) {$M\otimes_R R_{\mathcal{S}}$};
\node (tens t) at (0,0) {$M\otimes_R R_{\T}$};
\node (tors s) at (3,2) {$M_\mathcal{S}$};
\node (tors t) at (3,0) {$M_\T$};

\draw
(tens s) edge[->,>=angle 90] node[left] {$M\otimes_R \res^R_{\mathcal{S},\T}$} (tens t)
(tors s) edge[->,>=angle 90] node[right] {$\res^M_{\mathcal{S},\T}$} (tors t)
(tens s) edge[->,>=angle 90] node[above] {$\theta_{M,\mathcal{S}}$} (tors s)
(tens t) edge[->,>=angle 90] node[below] {$\theta_{M,\T}$} (tors t);

\end{tikzpicture}
\end{center}

Therefore $\theta_M$ gives a morphism of presheaves-on-a-basis from the presheaf underlying $\mathcal{M}_\otimes$ to the presheaf underlying $\mathcal{M}_\tors$.
\end{lemma}

\proof
We have a Yoneda isomorphism of $R$-modules $y_M:M\to(R_R,M)$ given by $y_M(m):r\mapsto mr$. Also, $Q_\T$ gives a morphism of abelian groups $(R_R,M)\to (Q_\T R,Q_\T M)=M_\T$. So define $\theta_{M,\T}(m\otimes \rho)= (Q_\T (y_M(m)))\circ \rho :Q_\T R\to Q_\T M$.

If $\mathcal{S}\subseteq\T$, take $m\otimes\rho$ in $M\otimes_R R_\mathcal{S}$. Following the diagram anticlockwise, $m\otimes \rho$ maps first to $m\otimes Q_{\mathcal{S},\T}(\rho)\in M\otimes_R R_\T$, and then to $(Q_\T (y_M(m)))\circ Q_{\mathcal{S},\T}(\rho)\in M_\T$. Going clockwise, $m\otimes \rho$ maps first to $(Q_{\mathcal{S}} (y_M(m)))\circ \rho\in M_\mathcal{S}$ and then to $Q_{\mathcal{S},\T}(Q_{\mathcal{S}} (y_M(m))\circ \rho)=(Q_\T (y_M(m)))\circ Q_{\mathcal{S},\T}(\rho)$, as before.\qed

Sheafifying, we therefore obtain a morphism of sheaves $\Theta_M:\mathcal{M}_\otimes\to\mathcal{M}_\tors$. This of course raises the question of what happens when we change modules along a map $f:M\to N$.

\begin{prop}
There is a natural transformation $\Theta$ from the tensor sheaf functor to the torsion sheaf functor, whose component at a module $M$ is $\Theta_M$.
\end{prop}

\proof
We must show that for any morphism $f:M\to N$, the following diagram commutes

\begin{center}
\begin{tikzpicture}

\node (tens m) at (0,2) {$\mathcal{M}_\otimes$};
\node (tors m) at (0,0) {$\mathcal{M}_\tors$};
\node (tens n) at (3,2) {$\mathcal{N}_\otimes$};
\node (tors n) at (3,0) {$\mathcal{N}_\tors$};

\draw
(tens m) edge[->,>=angle 90] node[left] {$\Theta_M$} (tors m)
(tens n) edge[->,>=angle 90] node[right] {$\Theta_N$} (tors n)
(tens m) edge[->,>=angle 90] node[above] {$f_\otimes$} (tens n)
(tors m) edge[->,>=angle 90] node[below] {$f_\tors$} (tors n);

\end{tikzpicture}
\end{center}

To do this, we show that for any basic open set $[A]$, the following diagram commutes

\begin{center}
\begin{tikzpicture}

\node (tens m) at (0,2) {$M\otimes_R R_{\T(A)}$};
\node (tors m) at (0,0) {$M_{\T(A)}$};
\node (tens n) at (5,2) {$N\otimes_R R_{\T(A)}$};
\node (tors n) at (5,0) {$N_{\T(A)}$};

\draw
(tens m) edge[->,>=angle 90] node[left] {$\theta_{M,\T(A)}$} (tors m)
(tens n) edge[->,>=angle 90] node[right] {$\theta_{N,\T(A)}$} (tors n)
(tens m) edge[->,>=angle 90] node[above] {$f\otimes_R R_{\T(A)}$} (tens n)
(tors m) edge[->,>=angle 90] node[below] {$f_{\T(A)}$} (tors n);

\end{tikzpicture}
\end{center}

Commutativity of this diagram establishes commutativity of the relevant diagram of pre\-sheaves-on-a-basis; as sheafification is functorial, it preserves commutativity of diagrams, and hence we obtain commutativity of the desired diagram of sheaves.

So to establish commutativity in our second diagram, we take $m\otimes \rho\in M\otimes_R R_{\T(A)}$. Following the diagram anticlockwise, we obtain first $(Q_{\T(A)} y_M(m))\circ \rho\in M_{\T(A)}$, and then $(Q_{\T(A)} f)\circ (Q_{\T(A)} y_M(m))\circ \rho\in N_{\T(A)}$. Going clockwise, $m\otimes \rho$ maps first to $f(m)\otimes \rho\in N\otimes_R R_{\T(A)}$, and then $(Q_{\T(A)} y_N(f(m)))\circ \rho\in N_{\T(A)}$. So we must show that $(Q_{\T(A)} f)\circ (Q_{\T(A)} y_M(m))=(Q_{\T(A)} y_N(f(m)))$.

Since $Q_{\T(A)}$ is functorial, $(Q_{\T(A)} f)\circ (Q_{\T(A)} y_M(m))=Q_{\T(A)}(f\circ y_M(m))$; so it suffices to prove that $f\circ y_M(m)=y_N(f(m))$. But this is precisely naturality of the Yoneda maps, completing the proof.\qed

\begin{cor}
For any torsion class $\T\in\Mod\mhyphen R$, there is a natural transformation $\theta_\T:-\otimes_R R_\T\to (-)_\T$.
\end{cor}

\proof
The component of $\theta_\T$ at the module $M$ is of course $\theta_{M,\T}$. The diagram whose commutativity needs checking is precisely the second diagram in the above proof.\qed

So we have two functors turning $R$-modules into sheaves of $\mathcal{O}_R$-modules, and a natural transformation between them. In the commutative noetherian case, we expect that torsion-theoretic localisation should be the same as localisation at a multiplicative set and hence that these two sheaf functors should coincide, so $\Theta$ should be an isomorphism. Indeed, we shall see a proof of this in Corollary \ref{commutative perfect}.

To address the question of when $\Theta$ is an isomorphism, we require the notion of a Gabriel filter. This is an alternative viewpoint on torsion-theoretic localisation, of which we give a brief overview based on Chapter VI of \cite{stenstrom}.\gap

Let $R$ be a ring and $\T$ a torsion class in $\Mod\mhyphen R$.  Then a module $M$ is $\T$-torsion if and only if every cyclic submodule of $M$ is $\T$-torsion. For, on the one hand, $\T$ is closed under subobjects, so any cyclic submodule of a $\T$-torsion module is $\T$-torsion; on the other hand, if $M$ is an $R$-module whose every cyclic submodule is $\T$-torsion, then $M$ can be expressed as a quotient of the direct sum of all its cyclic submodules, and so $M$ is $\T$-torsion.

Any cyclic module has the form $R/I$ for some right ideal $I$, so $\T$ is entirely determined by the set of right ideals $I$ such that $R/I$ is $\T$-torsion. We shall denote this set by $\mathcal{G}_\T$ and call it the {\bf Gabriel filter} associated to $\T$ (this terminology will be explained shortly). Recall that, for $I$ a right ideal and $r\in R$, $(I:r)$ denotes the right ideal $\{x\in R\mid rx\in I\}$; \textit{i.e.}, the annihilator of $r+I$ in the quotient module $R/I$.

\begin{lemma}[\cite{stenstrom}, \S\S VI.4, VI.5]
Let $R$ be any ring and $\T$ a torsion class in $\Mod\mhyphen R$. Let $\mathcal{G}_\T$ be the associated Gabriel filter:
\[\mathcal{G}_\T=\{I\leq R_R\mid R/I\in\mathcal{T}\}.\]
Then $\mathcal{G}_\T$ has the following three properties:
\begin{enumerate}
\item $\mathcal{G}_\T$ is a filter of right ideals of $R$ - \textit{i.e.}, it is closed under finite intersection and upwards inclusion;
\item If $I\in\mathcal{G}_\T$ and $r\in R$, then $(I:r)\in\mathcal{G}_\T$;
\item If $I$ is a right ideal and there is $J\in\mathcal{G}_\T$ such that for all $j\in J$, $(I:j)\in \mathcal{G}_\T$, then $I\in\mathcal{G}_\T$.
\end{enumerate}
\end{lemma}

Any collection of right ideals of $R$ satisfying the above 3 properties is called a {\bf Gabriel filter} on $R$, hence why $\mathcal{G}_\T$ is called the Gabriel filter associated to $\T$. Not only can we associate a Gabriel filter to any torsion theory on $\Mod\mhyphen R$, but we can also associate a torsion theory to any Gabriel filter $\mathcal{G}$ by declaring the cyclic torsion modules to be those of the form $R/I$ where $I\in\mathcal{G}$. We thus have the following:

\begin{thm}[\cite{stenstrom}, Theorem VI.5.1]
For any ring $R$ there is a bijective correspondence between torsion theories on $\Mod\mhyphen R$ and Gabriel filters on $R$.
\end{thm}

\begin{thm}[\cite{stenstrom}, Proposition XI.3.4]\label{perfect localisation}
Let $R$ be any ring and $\T$ a torsion class in $\Mod\mhyphen R$. Then the following are equivalent:
\begin{enumerate}
\item The functor $(\Mod\mhyphen R)/\T\to \Mod\mhyphen R_\T$ induced by passing $(-)_\T:M\mapsto M_\T$ to the quotient is an equivalence of categories;
\item The right adjoint inclusion $i_\T:(\Mod\mhyphen R)/\T\to \Mod\mhyphen R$ itself has a right adjoint;
\item The functor $(-)_\T:\Mod\mhyphen R\to \Mod\mhyphen R_\T:M\mapsto M_\T$ is exact and preserves coproducts;
\item The Gabriel filter $\mathcal{G}_\T$ has a filter base of finitely generated right ideals, and $(-)_\T$ is exact;
\item The natural transformation $\theta_\T:-\otimes_R R_\T\to (-)_\T$ is an isomorphism of functors;
\item For each $M\in\Mod\mhyphen R$, the kernel of the canonical map $M\to M\otimes_R R_\T$ is precisely $\tau_\T(M)$;
\item The restriction map $\res^R_\T:R\to R_\T$ is a ring epimorphism making $R_\T$ into a flat left $R$-module, and $\mathcal{G}_\T=\{I\leq R_R\mid \res^R_\T(I)R_\T=R_\T\}$.
\end{enumerate}
\end{thm}

We call a torsion class $\T$ satisfying the above equivalent conditions a {\bf perfect} torsion class.

\begin{thm}
Let $R$ be a right noetherian ring. Then the natural transformation $\Theta$ from the tensor sheaf functor to the torsion sheaf functor is a natural isomorphism if and only if every prime torsion class is perfect, if and only if for every prime torsion class $\T$ the functor $(-)_\T:\Mod\mhyphen R\to \Mod\mhyphen R_\T$ is exact.
\end{thm}

\proof
Since $\Theta$ is an isomorphism if and only if every component $\Theta_M$ is an isomorphism, it suffices to consider when $\Theta_M:\mathcal{M}_\otimes\to\mathcal{M}_\tors$ is an isomorphism of sheaves. A map of sheaves is an isomorphism if and only if the induced maps on stalks are all isomorphisms. The stalks are the localisations at torsionfree classes cogenerated by single indecomposable injectives; \textit{i.e.}, at prime torsion theories. So we see that $\Theta$ is an isomorphism if and only if $\theta_{M,\T}$ is an isomorphism for each module $M$ and prime torsion theory $\T$.

But $\theta_{M,\T}$ is precisely the component at $M$ of the natural transformation $\theta_\T: -\otimes_R R_\T\to (-)_\T$; so $\Theta$ is an isomorphism if and only if for each prime torsion theory $\T$, $\theta_\T$ is an isomorphism. By condition (5) of Theorem \ref{perfect localisation}, this occurs if and only if each prime torsion class is perfect.

Finally, we apply condition (4) of Theorem \ref{perfect localisation}. Since $R$ is right noetherian, every Gabriel filter has a filter base of finitely generated right ideals, so we see that $\Theta$ is an isomorphism if and only if $(-)_\T$ is exact for all prime torsion theories $\T$.\qed

We claimed above that over a commutative noetherian ring, the two sheaf functors are isomorphic along $\Theta$. We are now almost in a position to prove this, by proving that over such a ring all prime torsion classes are perfect; in fact, the stronger result holds that all torsion classes are perfect. First, though, we require some well-known preliminaries about torsion theories over commutative noetherian rings.

\begin{lemma}[\cite{gabriel}, Proposition V.5.10]\label{comm stable}
Let $R$ be a commutative noetherian ring. Then for any torsion theory $\T$ in $\Mod\mhyphen R$ and any prime ideal $p$, either $R/p\in\T$ or $R/p\in\mathcal{F}_\T$. Moreover, $R/p\in\T$ if and only if $E(R/p)\in \T$.
\end{lemma}

If a torsion theory $(\mathcal{T},\mathcal{F})$ has the property that every indecomposable injective is either in $\T$ or in $\mathcal{F}$, we say that it is a {\bf stable} torsion theory. The above Lemma shows that over a commutative noetherian ring, all torsion theories are stable.

The following result is well-known.

\begin{lemma}[\cite{stenstrom}, Example 2, of Section IX.1]\label{comm localisation is at a mult set}
Let $R$ be a commutative noetherian ring and $\T$ a torsion class in $\Mod\mhyphen R$. Then there is a multiplicative set $D\subseteq R$ such that $\T=\T_D$, the torsion class defined by
\[\T_D=\{M\in\Mod\mhyphen R\mid \forall m\in M\exists d\in D (md=0)\}.\]
\end{lemma}

Finally we are able to prove that for $R$ commutative noetherian, the natural transformation $\Theta$ is always an isomorphism.

\begin{cor}\label{commutative perfect}
Let $R$ be a commutative noetherian ring. Then every torsion class in $\Mod\mhyphen R$ is perfect.
\end{cor}

\proof
Let $\T$ be a torsion class in $\Mod\mhyphen R$. Then $\T=\T_D$ for some multiplicative set $D$, by Lemma \ref{comm localisation is at a mult set}, and the classical localisation at $D$ is an exact, full, and dense functor $\Mod\mhyphen R\to\Mod\mhyphen D^{-1}R$ with kernel exactly $\T$, so is equivalent to the torsion-theoretic localisation functor $Q_\T$, by the universal property of localisation (Proposition \ref{localisation}). More precisely, there is an equivalence of categories $F:(\Mod\mhyphen R)/\T\to \Mod\mhyphen D^{-1}R$ such that $F\circ Q_\T$ is the classical localisation functor.

This equivalence makes the adjoint inclusion $i_\T$ into the restriction of scalars functor $\Mod\mhyphen D^{-1}R\to\Mod\mhyphen R$, which has a right adjoint, namely the coinduced module functor $(D^{-1}R_R,-)$. So we meet condition (2) of Theorem \ref{perfect localisation}, and so $\T$ is perfect.\qed

So for a commutative noetherian ring, the two sheaves associated to a module coincide, and hence the two functors $\Mod\mhyphen R\to \sh(R)$ coincide too. Of course, these are simply the usual way of turning a module over a commutative ring into a sheaf over $\spec(R)$. We now turn to the consideration of noncommutative rings where these two sheaf functors coincide, making use of results from \cite{stenstrom}.

\begin{lemma}[\cite{stenstrom}, Proposition XI.3.3]
Let $R$ be any ring and $\mathcal{G}$ a Gabriel filter on $R$ having a filter base of projective right ideals. Let $\T$ be the torsion class associated to $\mathcal{G}$. Then $(-)_\T:\Mod\mhyphen R\to \Mod\mhyphen R_\T$ is exact.
\end{lemma}

\begin{cor}[\cite{stenstrom}, Corollary XI.3.6]
Let $R$ be a right noetherian, right hereditary ring. Then every torsion class in $\Mod\mhyphen R$ is perfect.
\end{cor}

\proof
Immediate from the above Lemma with part (4) of Theorem \ref{perfect localisation}.\qed

Therefore, for any right noetherian, right hereditary ring, such as a principal right ideal ring or the first Weyl algebra over a field of characteristic 0, the torsion sheaf functor and the tensor sheaf functor are naturally isomorphic. There is therefore a single sensible notion of the sheaf associated to a module, opening the way to exploration of further analogues with commutative algebraic geometry.

For a general ring, however, the localisations involved in the sheaf of finite-type localisations might fail to be perfect, in which case it is not clear which is the ``correct'' notion of the sheaf associated to a module. It may of course be that different contexts require considering either tensor sheaves or torsion sheaves.\gap

Recall that, given a ringed space $(X,\mathcal{O}_X)$, a sheaf of $\mathcal{O}_X$-modules $\mathcal{M}$ is {\bf quasicoherent} if it has everywhere a local presentation. That is, if for any point $x\in X$, there is a neighbourhood $U$, sets $I,J$, and an exact sequence of sheaves:
\[ \mathcal{O}_X^{(I)}|_U\to\mathcal{O}_X^{(J)}|_U\to \mathcal{M}|_U\to 0,\]
where $-|_U$ denotes the restriction of a sheaf on $X$ to a sheaf on $U$. Write $\qcoh(R)$ for the full subcategory of $\sh(R)$ consisting of the quasicoherent sheaves on $\injspec(R)$.

\begin{lemma}
Let $R$ be a ring such that each prime torsion class is perfect. Then the torsion sheaf functor (equivalently the tensor sheaf functor) $\Mod\mhyphen R\to \sh(R)$ lands in $\qcoh(R)$.
\end{lemma}

\proof
We will show that for any module $M$ there is in fact a global presentation for $\mathcal{M}_\tors$. Take a presentation for $M$ as an $R$-module:
\[ R^{(I)}\to R^{(J)}\to M\to 0.\]

Applying the torsion sheaf functor, we obtain a sequence of sheaves
\[\mathcal{O}_R^{(I)}\to \mathcal{O}_R^{(J)}\to \mathcal{M}_\tors\to 0;\]
we need only show that this sequence is exact. For this it suffices to show exactness on stalks. A stalk is given by localisation at a torsionfree class cogenerated by a single indecomposable injective; \textit{i.e.}, at a prime torsion theory, by Theorem \ref{torspec = injspec}. But, by hypothesis, these torsion theories are perfect, and so by part (3) of Theorem \ref{perfect localisation}, the localisation is exact.\qed

Therefore, for rings over which all prime torsion classes are perfect, we have a functor $\Mod\mhyphen R\to\qcoh(R)$. In the commutative case, this is an equivalence of categories. This result can certainly fail in the noncommutative case, as we now show.

\begin{ex}
Let $R=kA_2$, the path algebra over a field $k$ of the quiver $A_2$. Then the tensor sheaf functor is not an equivalence of categories.
\end{ex}

Recall Example \ref{kA2 counterexample}, where we showed that the ring of global sections of $\mathcal{O}_R$ was $k\oplus M_2(k)$. We show that $\qcoh(R)\cong \Mod\mhyphen (k\oplus M_2(k))$; \textit{i.e.}, that quasicoherent sheaves are equivalent to modules over the ring of global sections; since $k\oplus M_2(k)$ is not Morita equivalent to $R$, this proves that the tensor sheaf functor cannot be an equivalence.

First observe that, as $\injspec(R)$ is a 2-point discrete space, all sheaves are quasicoherent. Indeed, take any sheaf $\mathcal{M}\in\sh(R)$ and any point $E\in\injspec(R)$; then $\{E\}$ is open, and $\mathcal{M}_{\{E\}}$ is simply an $\mathcal{O}_R(\{E\})$-module, hence has a presentation. So $\qcoh(R)=\sh(R)$.

Write $E_1=(k\to 0)$ and $E_2=(k\to k)$ for the two indecomposable injective $R$-modules. Given a $(k\oplus M_2(k))$-module $M$, which can be naturally written as $M_1\oplus M_2$, for $M_1\in\Mod\mhyphen k$, $M_2\in\Mod\mhyphen M_2(k)$,, define a sheaf $\mathcal{M}$ by $\mathcal{M}(\{E_i\})=M_i$. A map $M\to N$ of $(k\oplus M_2(k))$-modules can be expressed as a pair of maps $(f_1,f_2)$, with $f_1:M_1\to N_1$, $f_2:M_2\to N_2$; this gives a morphism of sheaves $\mathcal{M}\to\mathcal{N}$. This defines a functor $\Mod\mhyphen(k\oplus M_2(k))\to\sh(R)$.

It is trivial to verify that this functor is quasi-inverse to the global sections functor, giving the desired equivalence of categories $\qcoh(R)\cong\Mod\mhyphen (k\oplus M_2(k))$.\qed

Although this shows that the global sections functor is not generally quasi-inverse to the tensor sheaf functor, we do at least have an adjunction between them, as we shall now show. For $M\in\Mod\mhyphen R$, let $\uline{M}$ denote the constant presheaf associated to $M$. Thus, $\uline{M}(U)=M$ for any open set $U$, and all restriction maps are the identity on $M$.

\begin{prop}
Let $R$ be a ring and let $\Gamma:\sh(R)\to\Mod\mhyphen R$ denote the global sections functor. Then the tensor sheaf functor is left adjoint to $\Gamma$.
\end{prop}

\proof
Since sheafification is left adjoint to the forgetful functor from sheaves to presheaves, it suffices to work with the presheaf $M\otimes_R \mathcal{O}_R$ assigning to an open set $U$ the $\mathcal{O}_R(U)$-module $M\otimes_R \mathcal{O}_R(U)$. For if $\mathcal{N}$ is any sheaf on $\injspec(R)$, then there is a natural isomorphism $(\mathcal{M}_\otimes,\mathcal{N})\cong (M\otimes_R\mathcal{O}_R,\mathcal{N})$; so we need only show the existence of a natural isomorphism $(M\otimes_R\mathcal{O}_R,\mathcal{N})\cong (M,\Gamma(\mathcal{N}))$.

A map $f:M\otimes_R \mathcal{O}_R\to \mathcal{N}$ consists of a map $f_U:M\otimes_R \mathcal{O}_R(U)\to \mathcal{N}(U)$ for each open set $U$, such that whenever $U\subseteq V$ the diagram below commutes.

\begin{center}
\begin{tikzpicture}

\node (mv) at (0,2) {$M\otimes_R \mathcal{O}_R(V)$};
\node (mu) at (0,0) {$M\otimes_R \mathcal{O}_R(U)$};
\node (nv) at (4,2) {$\mathcal{N}(V)$};
\node (nu) at (4,0) {$\mathcal{N}(U)$};

\draw
(mv) edge[->,>=angle 90] node[left] {$M\otimes_R \res^{\mathcal{O}_R}_{V,U}$} (mu)
(nv) edge[->,>=angle 90] node[right] {$\res^\mathcal{N}_{V,U}$} (nu)
(mv) edge[->,>=angle 90] node[above] {$f_V$} (nv)
(mu) edge[->,>=angle 90] node[below] {$f_U$} (nu);

\end{tikzpicture}
\end{center}

By tensor-hom adjunction, $f_V:M\otimes_R\mathcal{O}_R(V)\to\mathcal{N}(V)$ corresponds to the map $M\to (\mathcal{O}_R(V),\mathcal{N}(V))_R:m\mapsto f_V(m\otimes -)$, and similarly for $f_U$. Of course, $(\mathcal{O}_R(V),\mathcal{N}(V))$ is naturally isomorphic to $\mathcal{N}(V)$, and under this isomorphism, $f_V(m\otimes -)$ is identified with $f_V(m\otimes 1)$, which we shall denote $\hat{f}_U(m)$.

By commutativity of the above diagram, we see that $\res^\mathcal{N}_{V,U}\circ \hat{f}_V = \hat{f}_U\circ \res^{\mathcal{O}_R}_{V,U}$, which is the map sending $m\in M$ to $f_U(m\otimes\res^{\mathcal{O}_R}_{V,U}(1))=f_U(m\otimes 1)=\hat{f}(U)$.

So $f:M\otimes_R \mathcal{O}_R\to\mathcal{N}$ corresponds to a map $\hat{f}_U:M\to\mathcal{N}(U)$ for each open set $U$ such that whenever $U\subseteq V$ we have $\res^\mathcal{N}_{V,U}\hat{f}_V=\hat{f}_U$. But this is precisely the same as a map from the constant presheaf $\uline{M}$ to $\mathcal{N}$.

So sheaf maps $M\otimes_R\mathcal{O}_R\to\mathcal{N}$ correspond to presheaf maps $\uline{M}\to\mathcal{N}$; but these correspond naturally to maps $M\to\mathcal{N}(\injspec(R))$, since every component $\hat{f}_U$ of $\hat{f}:\uline{M}\to\mathcal{N}$ is just obtained from $\hat{f}_{\injspec(R)}$ by the formula $\hat{f}_U(m)=\res^\mathcal{N}_{\injspec(R),U}\circ \hat{f}_{\injspec(R)}$.

This establishes the isomorphism $(\mathcal{M}_\otimes,\mathcal{N})=(M,\Gamma(\mathcal{N}))$. Naturality follows from naturality of all the intermediate steps.\qed

\section{The Torsion Spectrum}\label{torsion}

Golan \cite{golan} discusses a number of topologies on the lattice of hereditary torsion theories in the module category over a noncommutative ring $R$ and a particular subset thereof, consisting of the prime torsion theories. This allows the definition of the `torsion spectrum' of a ring, which turns out, for $R$ noetherian, to be homeomorphic to the injective spectrum.

In fact, Golan's definitions, with a slight modification, work in an arbitrary Grothendieck category; so we work in this generality.

\subsection{Golan's Torsion Spectrum}

We begin by explaining the ideas of Golan \cite{golan}; the notation and terminology is significantly changed from that paper to fit in better with the other concepts in this paper. Fix a Grothendieck category $\mathcal{A}$.

An object $A$ of $\mathcal{A}$ is called {\bf torsion-critical} if every proper quotient of $A$ is $\mathcal{F}(A)$-torsion.

\begin{lemma}\label{torsion-critical}
Let $A$ be a torsion-critical object. Then
\begin{enumerate}
\item $A$ is uniform;
\item Any non-zero subobject of $A$ is torsion-critical;
\item For any non-zero subobject $B$ of $A$, $\mathcal{F}(B)=\mathcal{F}(A)$.
\end{enumerate}
\end{lemma}

\proof
\begin{enumerate}
\item Suppose for a contradiction that $B$ and $C$ are non-zero subobjects of $A$ with $B\cap C=0$. Then $A$ embeds in $A/B\oplus A/C$; but $A/B$ and $A/C$ are both $\mathcal{F}(A)$-torsion, hence so is $A$, a contradiction.
\item Let $B\leq A$ be non-zero. Note that, since $B\in\mathcal{F}(A)$, $\mathcal{F}(B)\subseteq\mathcal{F}(A)$. Suppose for a contradiction that $B$ has a proper quotient $C$ which is not $\mathcal{F}(B)$-torsion. Then
\[D:=\frac{C}{\tau_{\mathcal{F}(B)}(C)}\]
is a proper, non-zero quotient of $B$ which is $\mathcal{F}(B)$-torsionfree and hence $\mathcal{F}(A)$-torsion\-free. But then $D$ has the form $A'/A''$ for some $A''<A\leq M$, so $A/A''$ has a non-zero, $\mathcal{F}(A)$-torsionfree submodule, so is not $\mathcal{F}(A)$-torsion, a contradiction.
\item We prove the more general result that if $C$ is an essential subobject of an arbitrary object $B$, then $\mathcal{F}(C)=\mathcal{F}(B)$; by part (1.), this suffices. Since $C$ is essential in $B$, $E(C)=E(B)$, and the result then follows by Lemma \ref{cogenerating}, which says that $\mathcal{F}(C)=\mathcal{F}(E(C))$, and similarly for $B$.\qed
\end{enumerate}

The torsion theories of the form $\mathcal{F}(A)$ for $A$ torsion-critical are called {\bf prime torsion theories}. The set of all such is called the (right) {\bf torsion spectrum} of $\mathcal{A}$ and denoted $\torspec(\mathcal{A})$. As with the injective spectrum, when $\mathcal{A}=\Mod\mhyphen R$ for a ring $R$, we abuse notation to write $\torspec(R)=\torspec(\Mod\mhyphen R)$.

Golan's definition actually only considers those torsion theories of the form $\mathcal{F}(M)$ for $M$ a torsion-critical, cyclic $R$-module. However, by Lemma \ref{torsion-critical}, this is equivalent to our definition. For any non-zero, cyclic submodule of a torsion-critical module $M$ is also torsion-critical and cogenerates the same torsionfree class.

The torsion spectrum is endowed with a topology as follows. For $\mathcal{T}$ any torsion class, let $[\mathcal{T}]$ denote the set of prime torsion theories for which $\mathcal{T}$ is contained in the torsion class - \textit{i.e.}, the intersection of $\torspec(\mathcal{A})$ with the principal filter generated by $\mathcal{T}$ in the lattice of torsion classes. The set of all $[\mathcal{T}(A)]$ where $A$ ranges over $\mathcal{A}^\fp$ is a basis of open sets for a topology on $\torspec(\mathcal{A})$, called the finitary order topology by Golan. Henceforth, by $\torspec(\mathcal{A})$ we shall mean the set endowed with this particular topology. We denote by $(\mathcal{T})$ the complement of $[\mathcal{T}]$ in $\torspec(\mathcal{A})$.

\subsection{Torsion Theories and the Injective Spectrum}\label{torsion and injspec}

We now relate Golan's torsion spectrum to the injective spectrum.

\begin{lemma}\label{prime is ind cogen}
Let $\mathcal{A}$ be a locally noetherian Grothendieck category. Then a torsion theory $(\mathcal{T},\mathcal{F})$ is prime if and only if there is an indecomposable injective $E$ such that $\mathcal{F}=\mathcal{F}(E)$.
\end{lemma}

\proof
First we let $\mathcal{F}$ be a prime torsionfree class and show that it is cogenerated by a single indecomposable injective object. By definition, there is a torsion-critical object $A$ such that $\mathcal{F}=\mathcal{F}(A)$. By Lemma \ref{torsion-critical}, $A$ is uniform, so $E(A)$ is indecomposable. By Lemma \ref{cogenerating}, $\mathcal{F}(A)=\mathcal{F}(E(A))$, so indeed $\mathcal{F}$ has an indecomposable injective cogenerator.

Now let $E$ be indecomposable injective; we show that $\mathcal{F}(E)$ is prime. Since $\mathcal{A}$ is locally noetherian, $E$ has a non-zero noetherian subobject, $A$. By Lemma \ref{crit subobj}, $A$ has a non-zero subobject $B$ such that whenever $0<C\leq B$, $(B/C,E(B))=0$. But, by Lemma \ref{cogenerating}, this implies that $B/C\in\T_{\mathcal{F}(B)}$, so $B$ is torsion-critical. Since $E$ is uniform, $E=E(B)$, and so $\mathcal{F}(E)=\mathcal{F}(B)$ is a prime torsionfree class.\qed

\begin{thm}\label{torspec = injspec}
Let $\mathcal{A}$ be a locally noetherian Grothendieck category. Then there is a homeomorphism $\injspec(\mathcal{A})\to\torspec(\mathcal{A}):E\mapsto \mathcal{F}(E)$.
\end{thm}

\proof
We shall refer to this map as $h$ for the purposes of this proof. By Lemma \ref{prime is ind cogen}, $h$ is well-defined and surjective. To show injectivity, we must show that two indecomposable injectives which cogenerate the same torsionfree class are isomorphic; but this follows immediately from Lemma \ref{specialisation tfae} and Corollary \ref{t0}. For if $\mathcal{F}(E)=\mathcal{F}(F)$, then $E\leadsto F$ and $F\leadsto E$ in $\injspec(\mathcal{A})$, but $\injspec(\mathcal{A})$ is $T_0$, so $E\cong F$.

Now we show that $h$ is a homeomorphism. Since $h$ is bijective, it suffices to prove that $h([A])=[\T(A)]$ for any $A\in\mathcal{A}^\fp$. We have that $E\in [A]$ if and only if $(A,E)=0$, if and only if $A\in\mathcal{T}_{\mathcal{F}(E)}$, if and only if $\T(A)\subseteq\T_{\mathcal{F}(E)}$, if and only if $\T_{\mathcal{F}(E)}\in[\T(A)]$. But $h(E)$ is the torsion theory with torsion class $\T_{\mathcal{F}(E)}$, so this states precisely that $h(E)\in[\T(A)]$.\qed

Golan's definition of what he calls the finitary order topology on the torsion spectrum of a ring used the sets $[\mathcal{T}(M)]$ for $M$ cyclic as a basis, rather than for $M$ finitely presented. However, since the sets $[M]$ for $M$ cyclic form a basis for the topology on the injective spectrum, this shows that the $[\mathcal{T}(M)]$ for $M$ cyclic form a basis for the topology on the torsion spectrum (essentially by repeated applications of Lemma \ref{direct sum basic open sets}), so our definition is equivalent to Golan's. The approach adopted here, however, allowed us to work in the greater generality of an arbitrary Grothendieck category.

\subsection{Irreducibility and Sobriety}\label{irred and sober}

We now use this connection between the injective spectrum and torsion theories to develop further results about the topology on $\injspec(\mathcal{A})$. We begin with three well-known technical Lemmas.

\begin{lemma}\label{sum of serres}
Let $\{\mathcal{S}_i\mid i\in I\}$ be a collection of Serre subcategories of an abelian category. Then the Serre subcategory
\[\sum_{i\in I}\mathcal{S}_i\]
consists precisely of those objects admitting a finite filtration whose factors each lie in some $\mathcal{S}_i$.
\end{lemma}

\begin{lemma}[\cite{prest}, 11.1.14]\label{quotient of loc noeth}
Let $\mathcal{A}$ be a locally noetherian Grothendieck category. Then every torsion theory in $\mathcal{A}$ is of finite type and every torsion-theoretic quotient of $\mathcal{A}$ is also locally noetherian.
\end{lemma}

\begin{lemma}\label{ziegler closed = torsionfree}
Let $\mathcal{A}$ be a locally noetherian Grothendieck category. Then there is an in\-clusion-preserving bijection between Ziegler-closed subsets of the injective spectrum of $\mathcal{A}$ and hereditary torsionfree classes in $\mathcal{A}$.
\end{lemma}

\proof
Given $C\subseteq \injspec(\mathcal{A})$ Ziegler-closed, we associate the torsionfree class $\mathcal{F}(C)$. Given $\mathcal{F}$ a torsionfree class, we define $\mathcal{C}(\mathcal{F})=\injspec(\mathcal{A})\cap \mathcal{F}$, the set of indecomposable injectives in $\mathcal{F}$.

First we take $C\subseteq \injspec(\mathcal{A})$ Ziegler-closed and show that $\mathcal{C}(\mathcal{F}(C))=C$. It is clear that $C\subseteq \mathcal{C}(\mathcal{F}(C))$, so we prove the reverse inclusion. We have $C=\bigcap_{i\in I}\{A_i\}$ for some collection of finitely presented objects $A_i$, with $I$ some indexing set. If $E\in \mathcal{C}(\mathcal{F}(C))$, then $E\in \mathcal{F}(C)$, so $(T,E)=0$ for all $T\in\mathcal{T}_{\mathcal{F}(C)}$. Now, each $A_i$ has $(A_i,F)=0$ for all $F\in C$, so each $A_i\in\mathcal{T}_{\mathcal{F}(C)}$, so if $E\in \mathcal{C}(\mathcal{F}(C))$, then $(A_i,E)=0$ for all $i$ and so $E\in C$. Therefore $\mathcal{C}(\mathcal{F}(C))=C$, as required.

Now let $\mathcal{F}$ be a torsionfree class. We show first that $\mathcal{C}(\mathcal{F})$ is a Ziegler-closed set. Let $E\in \injspec(\mathcal{A})\take \mathcal{C}(\mathcal{F})$. Then $E$ is not torsionfree for $\mathcal{F}$, so there is an $\mathcal{F}$-torsion submodule $M$ of $E$ (which can be taken to be finitely presented, without loss of generality, as any non-zero subobject of $\tau_{\mathcal{F}}(E)$ will suffice for our argument). Then $(M,\mathcal{F})=0$, so for all $F\in \mathcal{C}(\mathcal{F})$, $(M,F)=0$, so $\mathcal{C}(\mathcal{F})\subseteq [M]$; but $(M,E)\neq 0$, so $E\in (M)\subseteq \injspec(\mathcal{A})\take \mathcal{C}(\mathcal{F})$, showing that $\injspec(\mathcal{A})\take \mathcal{C}(\mathcal{F})$ is Ziegler-open.

Now we show that $\mathcal{F}=\mathcal{F}(\mathcal{C}(\mathcal{F}))$. Since $\mathcal{C}(\mathcal{F})\subseteq\mathcal{F}$, certainly $\mathcal{F}(\mathcal{C}(\mathcal{F}))\subseteq\mathcal{F}$. Conversely, suppose $M\in\mathcal{F}$. Then $E(M)\in\mathcal{F}$, but all injectives are direct sums of indecomposable injectives, and $\mathcal{F}$ is closed under subobjects, so $E(M)$ is a direct sum of elements of $\mathcal{C}(\mathcal{F})$. Hence $E(M)\in\mathcal{F}(\mathcal{C}(\mathcal{F}))$, and hence so too is $M$.

Finally, it is clear that this preserves the inclusion ordering.\qed

Recall that a torsionfree class is prime if and only if it is cogenerated by a single indecomposable injective object (Lemma \ref{prime is ind cogen}).

\begin{cor}\label{sum of primes}
Every torsionfree class in a locally noetherian Grothendieck category $\mathcal{A}$ is a sum of prime torsionfree classes:
\[\mathcal{F}=\sum\{\mathcal{F}(E)\mid E\in\mathcal{F}\cap\injspec(\mathcal{A})\}.\]
\end{cor}

The following appears in \cite{prest}. The proof is not difficult, but requires a few extra technicalities, so we omit it.

\begin{lemma}[\cite{prest}, 11.1.10]\label{correspondence theorem for torsion classes}
Let $\mathcal{A}$ be a locally noetherian Grothendieck category and $\mathcal{T}$ a torsion class in $\mathcal{A}$. Then there is an inclusion preserving bijection between torsion classes in $\mathcal{A}$ which contain $\mathcal{T}$ and torsion classes in $\mathcal{A}/\mathcal{T}$.
\end{lemma}

Let us say that a torsion class is {\bf simple} if it properly contains no other torsion class except 0.

\begin{lemma}\label{simple torsion classes}
Let $\mathcal{A}$ be a locally noetherian Grothendieck category. Then for any simple object $S\in\mathcal{A}$, $\mathcal{T}(S)$ is a simple torsion class. Given two simple objects $S_1,S_2$, $\mathcal{T}(S_1)=\mathcal{T}(S_2)$ if and only if $S_1\cong S_2$.
\end{lemma}

\proof
The class $\mathcal{F}_{\mathcal{T}(S)}$ consists of those objects $F$ such that $(S,E(F))=0$; but $S$ is simple, so if $(S,E(F))\neq0$, then $S$ embeds in $E(F)$. Since $F$ is essential in $E(F)$, the image of $S$ in $E(F)$ has non-zero intersection with $F$, so is contained in $F$, by simplicity again. So $(S,E(F))\neq 0$ if and only if $(S,F)\neq 0$. So $\mathcal{F}_{\mathcal{T}(S)}$ consists of those $F$ such that $(S,F)=0$.

Therefore $\mathcal{T}(S)$ consists of those objects $T$ such that $(T,F)=0$, whenever $(S,F)=0$. Since $S$ is simple, $(S,F)=0$ precisely when $F$ does not contain $S$ as a subobject. So if any quotient of $T$ fails to contain $S$ as a submodule, that quotient is torsionfree but receives a map from $T$, a contradiction. So $\mathcal{T}(S)$ consists of objects whose every non-zero quotient has $S$ as a submodule.

But then any torsion class containing any non-zero object of $\mathcal{T}(S)$, being closed under subquotients, must contain $S$, and so contains all of $\mathcal{T}(S)$. So $\mathcal{T}(S)$ is a simple torsion class.

Now suppose that $S_1$ and $S_2$ are simple objects with $\mathcal{T}(S_1)=\mathcal{T}(S_2)$. Then $S_1\in\mathcal{T}(S_2)$, so $S_1$ has $S_2$ as a subobject, by the above. But $S_1$ is simple, so $S_1\cong S_2$, as claimed.\qed

\begin{thm}
Let $\mathcal{A}$ be a locally noetherian Grothendieck category and $(\mathcal{T},\mathcal{F})$ a torsion theory. Then the following are equivalent:
\begin{enumerate}
\item $(\mathcal{T},\mathcal{F})$ is prime;
\item $\mathcal{F}$ is $+$-irreducible in the lattice of torsionfree classes;
\item $\T$ is $\cap$-irreducible in the lattice of torsion classes.
\end{enumerate}
\end{thm}

\proof
($2.\Leftrightarrow 3.$): Since the lattice of torsion classes is dual to the lattice of torsionfree classes, this is obvious.\gap

($1.\Rightarrow 2.$): Let $\F$ be a prime torsionfree class, cogenerated by the indecomposable injective $E$. Suppose that $\mathcal{F}=\mathcal{F}_1+\mathcal{F}_2$ is the join of some torsionfree classes $\mathcal{F}_1$, $\mathcal{F}_2$ (necessarily contained in $\mathcal{F}$). Let $E_1, E_2$ be injective cogenerators for $\mathcal{F}_1, \mathcal{F}_2$ respectively. Then $E_1\times E_2$ cogenerates $\mathcal{F}$.

Since $E\in\mathcal{F}$, there is some cardinal $\lambda$ such that $E$ embeds in $(E_1\times E_2)^\lambda\cong E_1^\lambda\times E_2^\lambda$. Let $K_1, K_2$ be the kernels of the component maps $E\to E_1^\lambda$, respectively $E\to E_2^\lambda$. Then $K_1\cap K_2$ is the kernel of the embedding $E\hookrightarrow E_1^\lambda\times E_2^\lambda$, hence is 0. But $E$ is uniform, so either $K_1=0$ or $K_2$=0. Therefore $E$ is cogenerated by either $E_1$ or $E_2$, and so $\mathcal{F}=\mathcal{F}(E)$ is contained in $\mathcal{F}_1$ or $\mathcal{F}_2$.\gap

($3.\Rightarrow 1.$): Let $\T$ be a $\cap$-irreducible torsion class, with associated torsionfree class $\mathcal{F}$. We show that $\mathcal{A}/\mathcal{T}$ contains a unique simple object $S$ and that $i_\mathcal{T}E(S)$ is an indecomposable injective cogenerator for $\mathcal{F}$, showing that $\mathcal{F}$ is prime.

First note that, by Lemma \ref{correspondence theorem for torsion classes}, $\mathcal{A}/\T$ has at most one simple torsion class, since the intersection of two simple classes is necessarily 0, but 0 is $\cap$-irreducible in $\mathcal{A}/\T$.

Now note that, since $\mathcal{A}$ is locally noetherian, it certainly contains a noetherian object, $N$ say. Then $N$ has a maximal proper subobject, and hence a simple quotient. So $\mathcal{A}$ contains a simple object, $S$. Then by Lemma \ref{simple torsion classes}, $\mathcal{T}(S)$ is simple and, since two non-isomorphic simple objects must generate different torsion classes, we see that $\mathcal{A}/\T$ has exactly one simple object, $S$.

Since the coproduct of the injective hulls of the simple objects form a cogenerating set for $\mathcal{A}/\T$ (by the fact that every object has a simple subquotient - see \cite{lam}(Theorem 19.8), though the proof there deals specifically with module categories), we see that $E(S)$ is an injective cogenerator for all of $\mathcal{A}/\T$.

Consider the quotient functor $Q_\T:\mathcal{A}\to\mathcal{A}/\T$ and its right adjoint inclusion $i_\T$. Since $i_\T$ is fully faithful, it preserves indecomposables, and since it has an exact left adjoint, it preserves injectives. So $i_\T(E(S))$ is an indecomposable injective object of $\mathcal{A}$. We show that this object cogenerates $\mathcal{F}$, and therefore that $\mathcal{F}$ is prime.

Every object $F$ of $\mathcal{F}$ embeds in its localisation $i_\T Q_\T(F)$. Since $E(S)$ cogenerates $\mathcal{A}/\T$, there is some cardinal $\lambda$ such that $Q_\T(F)$ embeds in $E(S)^\lambda$. Since $i_\T$ is a right adjoint, it is left exact and preserves products, so $i_\T Q_\T(F)\hookrightarrow i_\T(E(S))^\lambda$. So every object of $\mathcal{F}$ is cogenerated by $i_\T(E(S))$.\qed

This allows us to identify points of $\injspec(\mathcal{A})=\torspec(\mathcal{A})$ purely in the lattice of torsion classes, without any reference to the actual objects of $\mathcal{A}$. However, it does not yet let us give a description of the topology, since this is given in terms of torsion classes generated by finitely presented objects. We will shortly address this, but first, we extract a Corollary from the above Proposition.

\begin{cor}\label{ziegler sober}
Let $\mathcal{A}$ be a locally noetherian Grothendieck category. Then the injective spectrum $\injspec(\mathcal{A})$ is sober in its Ziegler topology.
\end{cor}

\proof
Let $C\subseteq \injspec(\mathcal{A})$ be an irreducible Ziegler-closed set and let $\mathcal{F}=\mathcal{F}(C)$ be the torsionfree class it cogenerates. Since the lattice of Ziegler-closed subsets of $\injspec(\mathcal{A})$ is isomorphic to the lattice of torsionfree classes, by Lemma \ref{ziegler closed = torsionfree}, we see that $\mathcal{F}$ is $+$-irreducible, hence prime. So there is an indecomposable injective $E$ which cogenerates $\mathcal{F}$.

Now if $M$ is a finitely presented object of $\mathcal{A}$ and $(M,E)=0$, so $E\in [M]$, then $M\in\mathcal{T}_\mathcal{F}$, so $(M,\mathcal{F})=0$, so for all $F\in C$, $(M,F)=0$; so $C\subseteq [M]$. Therefore $C$ is contained in the Ziegler-closure of $E$. For the reverse inclusion, note that $E\in\mathcal{F}(C)\cap\injspec(R)$, which is $C$, by  Lemma \ref{ziegler closed = torsionfree}; since $C$ is Ziegler-closed and contains $E$, it contains the Ziegler-closure of $E$.\qed

Unfortunately, we are more interested in sobriety of the injective spectrum in its Zariski topology for the purposes of this thesis. This is partly resolved by the following result

\begin{prop}[\cite{prest}, Proposition 14.2.6]
Let $\mathcal{A}$ be a locally noetherian Grothendieck category and $(A)$ a non-empty basic closed set in $\injspec(\mathcal{A})$ (for some $A\in\mathcal{A}^\fp$). If $(A)$ is irreducible, then there is a generic point of $(A)$.
\end{prop}

Prest's proof of this in \cite{prest} is essentially purely topological. An alternative proof is possible using torsion-theoretic methods, by taking the sum of all torsion classes not containing $\T(A)$ and showing that it is prime and its corresponding indecomposable injective is generic. Another method, more model-theoretic in nature, involves taking an ultraproduct of all the indecomposable injectives in $(A)$ and showing that an indecomposable summand of that is generic in $(A)$. However, none of these methods has been able to address the existence of a generic point for a non-basic irreducible closed set, so it remains an open question whether the Zariski topology on $\injspec(\mathcal{A})$ is sober.

\begin{qn}
Is the injective spectrum sober?
\end{qn}

Sobriety in the Ziegler topology, as proved above, will however be useful in Section \ref{spectral spaces}.

\begin{lemma}
The Ziegler-closed sets of $\torspec(\mathcal{A})$ are precisely those of the form $[\T]$ for any torsion class $\T$.
\end{lemma}

\proof
A basis of closed sets for the Ziegler topology is given by the $[\T(A)]$ for $A$ finitely presented; so each Ziegler-closed set has the form
\[\bigcap_{i\in I} [\T(A_i)]\]
for some finitely presented objects $A_i$. For any torsion class $\T$, $\T$ contains all $\T(A_i)$ if and only if $\T$ contains $\sum_{i\in I}\T(A_i)$; \textit{i.e.}, we have
\[\bigcap_{i\in I} [\T(A_i)]=\left[\sum_{i\in I}\T(A_i)\right].\]

So all Ziegler-closed sets have the desired form. On the other hand, since in a locally noetherian category all torsion theories are of finite type and hence determined by their finitely presented objects, any torison class $\T$ can be written as the sum of the torsion classes generated by the finitely presented objects of $\mathcal{T}$, and so
\[[\T]=\left[\sum_{A\in \T\cap\mathcal{A}^{\mathrm{fp}}}\T(A)\right]=\bigcap_{A\in\T\cap\mathcal{A}^{\mathrm{fp}}}[\T(A)],\]
showing that every set of this form is Ziegler-closed.\qed

\begin{prop}\label{compact classes}
The torsion classes $\T(M)$ for $M\in\mathcal{A}^{\mathrm{fp}}$ are precisely those which are compact elements of the lattice of torsion classes; \textit{i.e.}, those $\T$ such that if $\T$ is contained in a sum of a set of torsion classes, it is already contained in the sum of some finite subset.
\end{prop}

\proof
First suppose that $\T$ is compact in the lattice. Since $\mathcal{A}$ is assumed to be locally noetherian, each torsion class $\T$ is determined by the finitely presented objects it contains. So
\[\T=\sum_{A\in\T^{\mathrm{fp}}} \T(A).\]
Since $\T$ is compact, there are some finitely presented objects $A_1,\hdots,A_n\in\T$ such that
\[\T=\sum_{i=1}^n\T(A_i)=\T\left(\bigoplus_{i=1}^n A_i\right),\]
so $\T$ is generated by a single finitely presented object.

Conversely, suppose $M$ is finitely presented and the torsion classes $\T_i$ for $i$ in some indexing set $I$ are such that
\[\T(M)\subseteq \sum_{i\in I}\T_i.\]
Intersecting with finitely presented objects, we have an inclusion of Serre subcategories of $\mathcal{A}^{\mathrm{fp}}$
\[\T(M)\cap\mathcal{A}^{\mathrm{fp}}\subseteq \sum_{i\in I}(\T_i\cap\mathcal{A}^{\mathrm{fp}}).\]
In particular, $M$ is contained in the right-hand side. By Lemma \ref{sum of serres}, therefore, $M$ admits a finite filtration, each of whose factors lies in some $\T_i\cap\mathcal{A}^{\mathrm{fp}}$; so there are finitely many $\T_i$ whose sum already contains $M$ and hence $\T(M)$. So $\T(M)$ is compact.\qed

So we now have a description of the torsion spectrum of a locally noetherian category purely in terms of the lattice of torsion classes. The points are the $\cap$-irreducible elements of the lattice, while a basis of open sets for the Zariski topology is the set of $[\T]$ where $\T$ ranges over compact elements of the lattice. The Ziegler topology has closed sets precisely the $[\T]$ for any $\T$.

\section{Spectral Spaces and Noetherianity}\label{spectral and noetherian}

In this section, we consider the consequences of two additional assumptions on $\injspec(R)$: that it be noetherian, and that it be spectral. We specialise to the injective spectrum of a right noetherian ring, rather than a general locally noetherian Grothendieck category, because certain compactness results can fail in the greater generality.

\subsection{Spectral Spaces}\label{spectral spaces}

Recall that a topological space $X$ is {\bf spectral} if it is compact, $T_0$, sober, and the the family of compact open sets of $X$ is closed under finite intersection and forms a basis of open sets for $X$. By a Theorem of Hochster \cite{hochster}, a space is spectral if and only if it is homeomorphic to the Zariski spectrum of some commutative ring. Moreover, given a spectral space $X$, there is an alternative ``dual'' topology on the same underlying set as $X$, where the complements of the compact open sets of $X$ are taken to be a basis of open sets for the dual topology. This dual space is also spectral, and its dual is the original topology on $X$.

It is not known whether the injective spectrum of a ring is spectral (in either of its topologies); however, it is `close enough' to spectral in its Ziegler topology to allow the dual topology to be defined, albeit without all the usual results holding, and this dual topology is precisely the Zariski topology (see subsection \ref{injspec intro}).

\begin{lemma}\label{reverse specialisation}
The specialisation order for the Ziegler topology on $\injspec(R)$ is simply the reverse ordering of the specialisation order in the Zariski topology.
\end{lemma}

\proof
Let $E,F\in\injspec(R)$ be indecomposable injectives. Then $E$ Zariski-specialises to $F$ if and only if every basic Zariski-closed set containing $E$ contains $F$. This means that for all finitely presented modules $M$ we have $(M,E)\neq 0$ implies $(M,F)\neq 0$. This occurs if and only if for all finitely presented $M$ we have $(M,F)=0$ implies $(M,E)=0$, which is precisely the statement that every basic Ziegler-closed set containing $F$ contains $E$. That is, that $F$ Ziegler-specialises to $E$.\qed

For any right noetherian ring $R$, $\injspec(R)$ is compact in its Ziegler topology \cite[5.1.11 \& 5.1.23]{prest} (this is the main point where we need to be over a ring, not just a locally noetherian Grothendieck category). By Proposition \ref{t0}, $\injspec(R)$ is $T_0$ in the Zariski topology for $R$ right noetherian; by Lemma \ref{reverse specialisation}, this implies that it is $T_0$ in the Ziegler topology too. Recall from subsection \ref{injspec intro} that the sets $(M)$ for $M\in\mod\mhyphen R$ are a basis of compact open sets for the Ziegler topology on $\injspec(R)$. By Corollary \ref{ziegler sober}, $\injspec(R)$ is sober in its Ziegler topology for any right noetherian ring $R$.

Therefore the only condition of a spectral space that can fail for the Ziegler topology on the injective spectrum of a right noetherian ring is that the intersection of compact open sets be compact open. If this condition holds, \textit{i.e.}, if $\injspec(R)$ is spectral in its Ziegler topology, then the Zariski topology, being the Hochster dual, is also spectral. In particular, this would prove sobriety of the injective spectrum in its Zariski topology. At present, no examples are known where the intersection of compact Ziegler-open sets of $\injspec(R)$ fails to be compact. So it is possible that the injective spectrum of a right noetherian ring is always a spectral space.

Given Hochster's result that all spectral spaces occur as Zariski spectra of commutative rings, if the injective spectrum of a right noetherian ring is always spectral, it would mean that a failure of commutativity cannot give anything new topologically, and that spectra of noncommutative rings differ only from those of commutative rings in the structure sheaf.

\begin{qn}
For $R$ right noetherian, is $\injspec(R)$ a spectral space? If not, are there necessary and/or sufficient conditions on $R$ for $\injspec(R)$ to be spectral?
\end{qn}

\subsection{Isolating Closed Sets}\label{isolating closed sets}

In order to prove statements about closed sets in injective spectra, it may be useful to isolate them; \textit{i.e.}, given a Zariski-closed set $C$ in the injective spectrum of some ring, to construct a ring, or at least a Grothendieck category, whose injective spectrum is homeomorphic to $C$. We will show that this can be done for basic closed sets if $\injspec(R)$ is spectral in its Ziegler topology, and for arbitrary closed sets if $\injspec(R)$ is also noetherian in its Zariski topology.

Consider first the case of a commutative noetherian ring $R$, so that the injective spectrum is the usual Zariski spectrum. A general closed set in $\spec(R)$ is $\spec(R/I)$ for some $I\ideal R$. In the injective spectrum, this corresponds to those indecomposable injectives $E$ which are the hulls of modules of the form $R/p$ for $p\in\spec(R)$ with $I\subseteq p$. This corresponds precisely to the basic closed set $(R/I)$ in $\injspec(R)$. Moreover, $\Mod\mhyphen R/I$ is a full subcategory of $\Mod\mhyphen R$, consisting of those modules which are quotients of direct sums of copies of $R/I$ - \textit{i.e.}, it is the full subcategory generated by $R/I$.

This suggests, then, for a basic closed set $(M)$ in the injective spectrum of an arbitrary noetherian ring $R$, to take the full subcategory of $\Mod\mhyphen R$ generated by $M$, in the hopes that the injective spectrum of this category will be homeomorphic to $(M)$. There is a problem, however; in general, this subcategory need not have a well-defined injective spectrum; indeed, it might not even be abelian. Following Wisbauer \cite[\S15]{wisbauer}, we consider the full subcategory $\sigma[M]$ of $\Mod\mhyphen R$ \textit{subgenerated} by $M$; \textit{viz.}, that consisting of all subquotients of direct sums of copies of $M$. This is the smallest Grothendieck subcategory of $\Mod\mhyphen R$ containing $M$.

We record in the next Theorem some results from \cite{wisbauer} that will be useful. Recall first that a module $E$ is said to be $M$-injective if for any submodule $N\leq M$ and map $f:N\to E$, $f$ extends to a map $M\to E$. For $R$-modules $A$ and $B$, define the trace of $A$ in $B$ by
\[\tr(A,B):=\sum_{f\in(A,B)}f(A).\]
If $A$ is fixed, then $\tr(A,-)$ is functorial (acting by restriction and corestriction on morphisms).

\begin{thm}[\cite{wisbauer}, 15.1, 16.3, 16.8, 17.9]\label{wisbauer}
For $M$ any $R$-module, we have:
\begin{enumerate}
\item The module
\[G_M:=\bigoplus\{U\leq M^{(\aleph_0)}\mid U\mbox{ is finitely generated}\}\]
is a generator for $\sigma[M]$ and the trace functor $\tr(G_M,-)$ is right adjoint to the inclusion $i_M:\sigma[M]\to \Mod\mhyphen R$;
\item The injective objects of $\sigma[M]$ are precisely those $M$-injective $R$-modules which lie in $\sigma[M]$;
\item For $N\in\sigma[M]$, if $E_R(i_MN)$ is its injective hull in $\Mod\mhyphen R$ and $E_M(N)$ is its injective hull in $\sigma[M]$, then $E_M(N)=\tr(M,E_R(i_MN))$.
\end{enumerate}
\end{thm}

In particular, for $R$ right noetherian and $M$ finitely presented, the summands of $G_M$ are noetherian, so $\sigma[M]$ has a generating set of noetherian objects; \textit{i.e.}, $\sigma[M]$ is locally noetherian.

We now consider how to use this to isolate a closed set. We keep the notation introduced above.

\begin{thm}\label{sigma m bijection}
Let $M$ be a finitely presented $R$-module. Then there exists a bijection $j:\injspec(\sigma[M])\to (M)\subseteq\injspec(R):F\mapsto E_R(i_MF)$, with inverse $j^{-1}=\tr(M,-)$.
\end{thm}

\proof
Let $F\in (M)$; then there is a non-zero map $M\to F$, so $\tr(M,F)\neq 0$. Since $i_M\tr(M,F)$ is a non-zero submodule of $F$, $F=E_R(i_M\tr(M,F))$; note that this already proves that $j\circ j^{-1}$ is the identity on $(M)$, under the assumption that $j$ and $j^{-1}$ are well-defined. If $\tr(M,F)$ were decomposable, then, since $i_M$ is fully faithful, so too would be $i_M\tr(M,F)$, and so $E_R(i_M\tr(M,F))$ would be decomposable, which is a contradiction. So $\tr(M,F)$ is non-zero and indecomposable, and is injective by part (3) of Theorem \ref{wisbauer}. So $\tr(M,-)$ does indeed give a well-defined function $(M)\to \injspec(\sigma[M])$.

For any $F\in\injspec(\sigma[M])$, $F$ is uniform in $\sigma[M]$. Since $\sigma[M]$ is a full subcategory of $\Mod\mhyphen R$ and is closed under submodules, this implies that $i_MF$ is uniform in $\Mod\mhyphen R$, so $E_R(i_MF)$ is indecomposable. So $j$ is well-defined.

Finally, by part (3) of Theorem \ref{wisbauer}, we see that $\tr(M,-)\circ j$ is the identity function on $\injspec(\sigma[M])$.\qed

Henceforth, we identify $\sigma[M]$ as a category in its own right with its image under $i_M$, except where this might cause confusion. Note that some categorical constructions do depend on whether we are working in $\sigma[M]$ or in $\Mod\mhyphen R$; for instance, products in $\sigma[M]$ are not the same as in $\Mod\mhyphen R$ \cite[{}15.1]{wisbauer}. For $N\in\sigma[M]$, the notation $[N]_R$ will refer to the basic open set of $\injspec(R)$ determined by $N$, whereas $[N]_M$ will denote the basic open set in $\injspec(\sigma[M])$. A similar convention will be adopted for basic closed sets.\gap

We now wish to prove that $j$ is a homeomorphism. This is where we require $\injspec(R)$ to be spectral in the Ziegler topology. We first require the following

\begin{lemma}\label{tr(m,e)=tr(gm,e)}
Let $E$ be an indecomposable injective $R$-module and $M$ any $R$-module. Then $\tr(M,E)=\tr(G_M,E)$.
\end{lemma}

\proof
Since $M$ is a summand of $G_M$, we certainly have that $\tr(M,E)\subseteq \tr(G_M,E)$. Let $i:\tr(G_M,E)\to E$ be the inclusion; then, given $A\leq M$ and $f:A\to \tr(G_M,E)$, $i\circ f$ is a map $A\to E$ in $\Mod\mhyphen R$. This therefore extends to a map $g:M\to E$, whose image must lie in $\tr(M,E)\subseteq\tr(G_M,E)$; so the corestriction of $g$ is a morphism $M\to \tr(G_M,E)$ extending $f$. Therefore $\tr(G_M,E)$ is $M$-injective and hence, by part (2) of Theorem \ref{wisbauer}, $\tr(G_M,E)$ is an injective object of $\sigma[M]$.

Now, since $\tr(M,E)$ is injective in $\sigma[M]$, by part (3) of the same Theorem (or, indeed, by a minor adaptation of the proof just given for $\tr(G_M,E)$), it is a summand of $\tr(G_M,E)$. So to prove equality, it suffices to show that $\tr(G_M,E)$ is indecomposable. But, in $\Mod\mhyphen R$, it is a subobject of the uniform module $E$, so it certainly is indecomposable.\qed

\begin{thm}\label{isolate}
Suppose that $\injspec(R)$ is spectral in its Ziegler topology. Then the above bijection $j$ is a homeomorphism when $(M)$ has the subspace topology inherited from the Zariski topology on $\injspec(R)$.
\end{thm}

\proof
First we prove that $j^{-1}$ is continuous. By Theorem 25.1 of \cite{wisbauer}, finitely presented $R$-modules in $\sigma[M]$ are finitely presented as objects of $\sigma[M]$, and for any $N\in\sigma[M]^\mathrm{fp}$, $N$ is a subquotient of a finite direct sum of copies of $M$. But for $M$ finitely presented over a right noetherian ring $R$, this implies that $N$ is finitely presented as an $R$-module. So $\sigma[M]^\mathrm{fp}=\sigma[M]\cap \mod\mhyphen R$.

Now, a basic open set of $\injspec(\sigma[M])$ has the form $[N]_M$ for $N\in\sigma[M]^\mathrm{fp}$. We show that $j[N]_M=[N]_R\cap (M)$, which suffices. This amounts to showing that, for $F\in\injspec(\sigma[M])$, $(N,F)=0$ if and only if $(N,E_R(F))=0$.

Since $(N,-)$ is left-exact and $F\subseteq E_R(F)$, we certainly have that if $(N,E_R(F))=0$, then $(N,F)=0$. Conversely, if $f:N\to E_R(F)$ is non-zero, then $f(N)\cap F\neq 0$, since $E_R(F)$ is uniform, so there is $n\in N$ such that $0\neq f(n)\in F$. Then $f|_{nR}:nR\to F$ is a non-zero morphism in $\sigma[M]$, but $F$ is injective in $\sigma[M]$, so this extends to a non-zero morphism $N\to F$. Note that we have shown that $j^{-1}$ is continuous not only in the Zariski topology, but also in the Ziegler topology; this will be essential for proving that $j$ is Zariski-continuous.\gap

Now we prove continuity of $j$. Let $N\in\mod\mhyphen R$; we show that $j^{-1}((N)_R\cap (M)_R)$ is closed. Note that $j^{-1}((N)_R\cap (M)_R)=\{F\in\injspec(\sigma[M])\mid (N,E_R(F))\neq 0\}$.

First we deal with the case where $N\in\sigma[M]$, so $(N)_R\subseteq (M)_R$. Then we have $(N,E_R(F))=(N,Tr(G_M,E_R(F))$, by part (1) of Theorem \ref{wisbauer}, and this is $(N,F)$, by Theorem \ref{sigma m bijection} and Lemma \ref{tr(m,e)=tr(gm,e)}, so $j^{-1}(N)_R=(N)_M$ in this case.

Now we deal with general $N$. Let $S$ be the set of subquotients of $N$ which lie in $\sigma[M]$. We show first that
\begin{equation}
j^{-1}((N)_R\cap (M)_R)=j^{-1}\left(\bigcup_{L\in S}(L)_R\right)=\bigcup_{L\in S}j^{-1}(L)_R.\tag{$\star$}
\end{equation}

The second equality is a standard fact about images of unions. For the first equality, take $F\in \injspec(\sigma[M])$. Then $(N,E_R(F))\neq0$ if and only if there is a subquotient $L$ of $N$ which embeds in $F$. Since $\sigma[M]$ is closed under subobjects, $L\in\sigma[M]$, so $L\in S$, and so $F\in (L)_R$. Conversely, if any subquotient $L$ of $N$ has a non-zero map to $F$, then $(L,E_R(F))\neq 0$, and by injectivity we conclude $(N,E_R(F))\neq 0$. This proves the claim.

Now each $L\in S$ lies in $\sigma[M]$, so $j^{-1}(L)_R=(L)_M$, by what we showed above. Moreover, $(N)_R$ and $(M)_R$ are compact open in the Ziegler topology, so $(N)_R\cap (M)_R$ is compact, since the Ziegler topology is assumed to be spectral. Since $j^{-1}$ is Ziegler-continuous and the continuous image of a compact set is compact, $j^{-1}(N)_R$ is Ziegler-compact. But each $(L)_M$ is Ziegler-open in $\injspec(\sigma[M])$, so the union in equation $(\star)$ can be replaced by a finite union.

So $j^{-1}(N)_R$ is equal to a finite union of sets of the form $j^{-1}(L)_R=(L)_M$ for $L\in\sigma[M]^\fp$, hence is Zariski-closed.\qed

So if $\injspec(R)$ is spectral in its Ziegler topology, then for any basic closed set $(M)$ there is a locally noetherian category $\sigma[M]$ whose injective spectrum is homeomorphic to $(M)\subseteq\injspec(R)$. What about arbitrary closed sets? These are intersections of basic closed sets; if $\injspec(R)$ is noetherian (in its Zariski topology), then every closed set is a \textit{finite} intersection of basic closed sets, and is therefore a single basic closed set, by the spectrality assumption. So for $\injspec(R)$ Zariski-noetherian and Ziegler-spectral, the above result covers all closed sets. Note, however, that subsection 4.4 of \cite{gulliver1} exhibits a noetherian ring whose injective spectrum is not noetherian.

\end{document}